\documentclass[a4paper,11pt]{article}
\usepackage[latin1]{inputenc} 
\usepackage[T1]{fontenc}
\usepackage{lmodern}

\usepackage{amsthm,amsmath,amsfonts,stmaryrd,bbm,mathrsfs}
\usepackage{mathtools}
\usepackage{enumitem}
\usepackage[colorlinks=true, linkcolor=blue, citecolor=blue,pdfstartview=FitH]{hyperref}

\usepackage{graphicx,tikz,color}
\usetikzlibrary{shapes}
\usetikzlibrary{patterns}

\usepackage[top=4cm, bottom=4cm, left=3cm, right=3cm]{geometry}
\usepackage{titlesec} 

\usepackage[french,english]{babel}

\usepackage{caption}

\makeatletter
\renewcommand\theequation{\thesection.\arabic{equation}}
\@addtoreset{equation}{section}
\makeatother

\renewcommand{\captionfont}{\footnotesize}
\renewcommand{\captionlabelfont}{\footnotesize}

\setenumerate[1]{label=(\roman*), font = \normalfont}

\let\originalleft\left
\let\originalright\right
\renewcommand{\left}{\mathopen{}\mathclose\bgroup\originalleft}
\renewcommand{\right}{\aftergroup\egroup\originalright}

\newcommand{\N}{\mathbb{N}}
\newcommand{\Z}{\mathbb{Z}}
\newcommand{\Q}{\mathbb{Q}}
\newcommand{\R}{\mathbb{R}}
\newcommand{\C}{\mathbb{C}}
\newcommand{\T}{\mathbb{T}}
\renewcommand{\P}{\mathbb{P}}
\newcommand{\E}{\mathbb{E}}

\newcommand{\Ec}[1]{\mathbb{E} \left[#1\right]}
\newcommand{\Ep}[1]{\mathbb{E} \left(#1\right)}
\newcommand{\Pc}[1]{\mathbb{P} \left[#1\right]}
\newcommand{\Pp}[1]{\mathbb{P} \left(#1\right)}
\newcommand{\Ecsq}[2]{\mathbb{E} \left[#1\mathrel{}\middle|\mathrel{}#2\right]}
\newcommand{\Epsq}[2]{\mathbb{E} \left(#1\mathrel{}\middle|\mathrel{}#2\right)}
\newcommand{\Pcsq}[2]{\mathbb{P} \left[#1\mathrel{}\middle|\mathrel{}#2\right]}
\newcommand{\Ppsq}[2]{\mathbb{P} \left(#1\mathrel{}\middle|\mathrel{}#2\right)}
\newcommand{\tEc}[1]{\widetilde{\mathbb{E}} \left[#1\right]}
\newcommand{\tPp}[1]{\widetilde{\mathbb{P}} \left(#1\right)}
\newcommand{\tEcsq}[2]{\widetilde{\mathbb{E}} 
						\left[#1\mathrel{}\middle|\mathrel{}#2\right]}
\newcommand{\tPpsq}[2]{\widetilde{\mathbb{P}} 
						\left(#1\mathrel{}\middle|\mathrel{}#2\right)}
\newcommand{\Eci}[2]{\mathbb{E}_{#1} \left[#2\right]}
\newcommand{\Ppi}[2]{\mathbb{P}_{#1} \left(#2\right)}
\newcommand{\Ppsqi}[3]{\mathbb{P}_{#1} \left(#2\mathrel{}\middle|\mathrel{}#3\right)}
\newcommand{\Ecsqi}[3]{\mathbb{E}_{#1} \left[#2\mathrel{}\middle|\mathrel{}#3\right]}

\newcommand{\1}{\mathbbm{1}}

\newcommand{\cF}{\mathcal{F}}
\newcommand{\cN}{\mathcal{N}}
\newcommand{\cP}{\mathcal{P}}
\newcommand{\cD}{\mathcal{D}}
\newcommand{\cQ}{\mathcal{Q}}
\newcommand{\cC}{\mathcal{C}}
\newcommand{\cH}{\mathcal{H}}
\newcommand{\cE}{\mathcal{E}}
\newcommand{\cU}{\mathcal{U}}
\newcommand{\cM}{\mathcal{M}}
\newcommand{\cL}{\mathcal{L}}
\newcommand{\cT}{\mathcal{T}}
\newcommand{\cG}{\mathcal{G}}
\newcommand{\cJ}{\mathcal{J}}
\newcommand{\cS}{\mathcal{S}}

\newcommand{\sF}{\mathscr{F}}
\newcommand{\sN}{\mathscr{N}}
\newcommand{\sC}{\mathscr{C}}
\newcommand{\sU}{\mathscr{U}}
\newcommand{\sE}{\mathscr{E}}
\newcommand{\sG}{\mathscr{G}}
\newcommand{\sT}{\mathscr{T}}
\newcommand{\sM}{\mathscr{M}}

\newcommand{\e}{\mathrm{e}}
\newcommand{\diff}{\mathop{}\mathopen{}\mathrm{d}}
\DeclareMathOperator{\Var}{Var}

\newcommand{\abs}[1]{\left\lvert#1\right\rvert}
\newcommand{\norme}[1]{\left\lVert#1\right\rVert}

\newcommand{\petito}[1]{o\mathopen{}\left(#1\right)}
\newcommand{\grandO}[1]{O\mathopen{}\left(#1\right)}

\newcommand{\enstq}[2]{\left\{#1\mathrel{}\middle|\mathrel{}#2\right\}}
\newcommand{\restreinta}{\mathclose{}|\mathopen{}}

\title{Velocity of the $L$-branching Brownian motion}

\author{Michel \bsc{Pain}\footnote{LPMA, UPMC and DMA, ENS.}}

\theoremstyle{plain}  
\newtheorem{thm}{Theorem}[section]
\newtheorem{prop}[thm]{Proposition}
\newtheorem{lem}[thm]{Lemma}
\newtheorem{cor}[thm]{Corollary}

\theoremstyle{definition}

\theoremstyle{remark}
\newtheorem{rem}[thm]{Remark}

\begin{document}

\maketitle

\vspace{-0.5cm}

\begin{abstract}
\noindent We consider a branching-selection system of particles on the real line that evolves according to the following rules: each particle moves according to a Brownian motion during an exponential lifetime and then splits into two new particles and, when a particle is at a distance $L$ of the highest particle, it dies without splitting.
This model has been introduced by Brunet, Derrida, Mueller and Munier \cite{bdmm2006-1} in the physics literature
and is called the $L$-branching Brownian motion. 
We show that the position of the system grows linearly at a velocity $v_L$ almost surely and we compute the asymptotic behavior of $v_L$ as $L$ tends to infinity: 
\[v_L = \sqrt{2} - \frac{\pi^2}{2\sqrt{2}L^2} + \petito{\frac{1}{L^2}},\]
as conjectured in \cite{bdmm2006-1}.
The proof makes use of results by Berestycki, Berestycki and Schweinsberg \cite{bbs2013} concerning branching Brownian motion in a strip.
\end{abstract}

\vspace{0.3cm}

{\small 

}

\section{Introduction}

The branching Brownian motion (or BBM) is a branching Markov process whose study dates back to \cite{inw68-2}.
It has been the subject of a large literature, especially for its connection with the F-KPP equation, highlighted by McKean \cite{mckean75}.
It is defined as follows.
Initially, there is a single particle at the origin. 
Each particle moves according to a Brownian motion, during an exponentially distributed time and then splits into two new particles, which start the same process from their place of birth.
Every particle behaves independently of the others and the system goes on indefinitely.
We study here a branching Brownian motion with selection. 
A particle's position corresponds to its survival capacity and reproductive success (biologists call it fitness): 
it changes during the particle's life because of mutations and is then transmitted to the particle's children.
The selection tends to eliminate the lowest particles, that have a too small fitness value by comparison with the best ones.
Thus, we consider a system of particles evolving as before, but where in addition a particle dies as soon as it is at a distance $L$ of the highest particle alive at the same time.
This system is called the $L$-branching Brownian motion or $L$-BBM.

\subsection{Statement of the results}

\begin{figure}[ht]
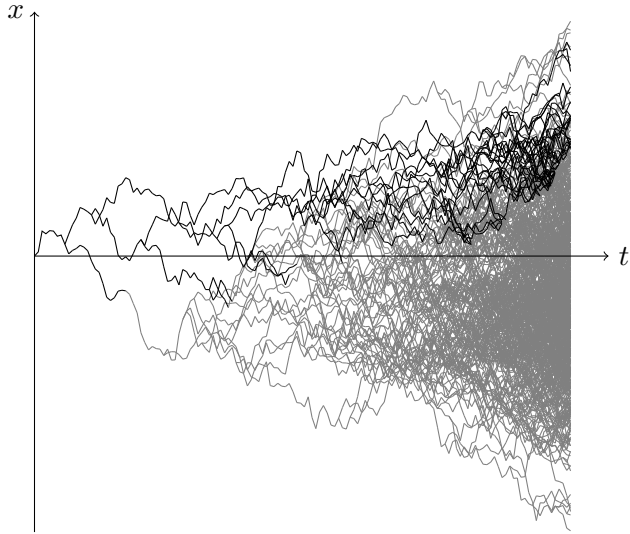
 
\centering
 
\caption{
Simulation of a BBM starting with a single particle at 0 and of the coupled $L$-BBM with $L=3$ between times 0 and 7. 
Particles in black belong to the $L$-BBM and particles in grey belong only to the BBM without selection and not to the $L$-BBM.
Note that quickly some particles killed by selection have descendants with a better fitness value than the particles of the $L$-BBM.} 
\label{figure L-BBM et BBM}
\end{figure}

First of all, we need to define the branching Brownian motion for a more general initial condition.
Let $\cC \coloneqq \bigcup_{n\in\N^*} \R^n$ be the set of configurations.
For each $\xi \in \cC$, $\xi = (\xi_1,\dots,\xi_n)$, we define the branching Brownian motion starting from this configuration as before, but with $n$ particles at time 0 positioned at $\xi_1,\dots,\xi_n$.
We denote by $M(t)$ the number of particles in the BBM at time $t$ and by $X_1(t), \dots, X_{M(t)}(t)$ their positions.
We will say that $(X_k(t), 1 \leq k \leq M(t))_{t\geq 0}$ is a branching Brownian motion although this notation does not contain the genealogy of the process.
We work on a measure space $(\Omega, \sF, (\P_\xi)_{\xi \in \cC})$ such that for each $\xi \in \cC$, under $\P_\xi$, $(X_k(t), 1 \leq k \leq M(t))_{t\geq 0}$ is a branching Brownian motion starting from the configuration $\xi$.
We equip this space with the canonical filtration associated to the branching Brownian motion, denoted by $(\sF_t)_{t\geq 0}$.
For each $t\geq0$, let $X(t) \coloneqq (X_1(t), \dots, X_{M(t)}(t)) \in \cC$ be the configuration of the particles of the BBM living at time $t$.
For $\xi = (\xi_1,\dots,\xi_n)$ and $\xi' = (\xi'_1,\dots,\xi'_m)$ configurations,
we say that $\xi' \subset \xi$ if $\sum_{i=1}^m \delta_{\xi'_i} \leq \sum_{i=1}^n \delta_{\xi_i}$ as measures
and we define $\max \xi \coloneqq \max_{1\leq i \leq n} \xi_i$.
Thus, $\max X(t)$ is the position of the highest particle of the BBM living at time $t$.
Then, we define the $L$-branching Brownian motion on the same space $(\Omega, \sF,  (\sF_t)_{t\geq 0}$, $(\P_\xi)_{\xi \in \cC})$, by coupling it canonically with the standard branching Brownian motion: for each realization of the BBM, we can define a realization of the $L$-BBM by killing every particle that is at a distance greater than $L$ from the highest particle of the $L$-BBM.
We denote by $M^L(t)$ the number of particles in the $L$-BBM at time $t$ and by $X^L_1(t), \dots, X^L_{M^L(t)}(t)$ their position.
Let $X^L(t) \coloneqq (X^L_1(t), \dots, X^L_{M^L(t)}(t))$.
By definition, we have the inclusion $X^L(t) \subset X(t)$ for all $t\geq 0$ and $\omega \in \Omega$ and,
for each $\xi \in \cC$, under $\P_\xi$, $(X^L_k(t), 1 \leq k \leq M^L(t))_{t\geq 0}$ is an $L$-branching Brownian motion starting from the configuration $\xi$.
See Figure \ref{figure L-BBM et BBM}.
Our aim here is to study the asymptotic behavior of the position of the highest particle of the $L$-branching Brownian motion at time $t$, that is $\max X^L(t)$.
For the standard branching Brownian motion, it is known \cite{chtww95} that 
\[
\frac{1}{t} \max X(t)
\underset{t\to\infty}{\longrightarrow}
\sqrt{2}
\]
almost surely.
So, the speed of the highest particle of the $L$-BBM, whenever it exists, has to be no more than $\sqrt{2}$, because there are less particles in the $L$-BBM than in the BBM owing to selection.
Moreover, it is clear that if this velocity exists, then all particles of the $L$-BBM moves at the same speed.
Our first result ensures that the velocity of the $L$-branching Brownian motion is well-defined and does not depend on the initial configuration.
The same result is proved in a similar way by Derrida and Shi in their article \cite{derridashi2016arxiv} in preparation at the time of writing of this paper, but, nonetheless, we present it here in order to work legitimately on $v_L$ in what follows.
\begin{prop} \label{proposition existence}
For each $L>0$, there exists $v_L \in \R$ such that for each $\xi \in \cC$ we have the following convergence
\[
\frac{1}{t} \max X^L(t)
\underset{t\to\infty}{\longrightarrow}
v_L
\]
$\P_\xi$-almost surely.
\end{prop}
We focus in this paper on the behavior of $v_L$ as $L$ tends to infinity.
It means that the selection effect vanishes, so we can expect that $v_L$ tends to $\sqrt{2}$, which is the asymptotic velocity of the highest particle of the BBM without selection.
Furthermore, we are interested in the asymptotic order of $\sqrt{2} - v_L$, which permits to estimate the slowdown due to selection.
The main result of this work shows this second term in the asymptotic expansion of $v_L$, validating a conjecture of Brunet, Derrida, Mueller and Munier \cite{bdmm2006-1} (see Subsection \ref{subsection motivations} for more details).
\begin{thm} \label{theorem}
We have the following asymptotic behavior:
\[ v_L = \sqrt{2} - \frac{\pi^2}{2\sqrt{2}L^2} + \petito{\frac{1}{L^2}},\]
as $L$ tends to infinity.
\end{thm}
This is analogous to the result of Bérard and Gouéré \cite{berardgouere2010} for the $N$-branching random walk, where the selection imposes a constant population size $N$.
One can expect that the population of the $L$-BBM is of order $\e^{c L}$ with $c$ some positive constant and, then, the result of Bérard and Gouéré suggests by taking $N = \e^{c L}$ that the first correction term for $v_L$ must be of order $1/L^2$.
Actually, we will see in the next subsection that Brunet, Derrida, Mueller and Munier \cite{bdmm2006-1} conjecture that $c = \sqrt{2}$ and it leads exactly to the term $- \pi^2/2\sqrt{2}L^2$.
However, the strategy to prove Theorem \ref{theorem} will neither be to use comparisons with the $N$-BBM nor to control precisely the population size of the $L$-BBM.

\subsection{Motivations} \label{subsection motivations}

The $L$-branching Brownian motion has been introduced by Brunet, Derrida, Mueller and Munier \cite{bdmm2006-1} in order to describe the effect of a white noise on the F-KPP equation.
The F-KPP equation (for Fisher \cite{fisher37} and Kolmogorov, Petrovsky, Piscounov \cite{kpp37}) 
\begin{equation} \label{equation F-KPP}
\partial_t h 
= \frac{1}{2} \partial^2_x h + h(1-h)
\end{equation}
is a traveling wave equation that describes how a stable phase ($h=1$) invades an unstable phase ($h=0$). 
Depending on the initial condition, the front between the two phases can travel at any velocity $v$ larger than a minimal velocity $v_\text{min} = \sqrt{2}$ (which is as well the asymptotic velocity of the highest particle of standard BBM).
This equation often represents a large-scale limit or a mean-field description of some microscopic discrete stochastic processes.
In order to understand the fluctuations that appear at the microscopic scale with a finite number of particles, 
one might consider instead the F-KPP equation with noise \cite{brunetderrida2001,bdmm2006-1,bdmm2006-2}
\begin{equation} \label{equation F-KPP with noise}
\partial_t h 
= \frac{1}{2} \partial^2_x h + h(1-h) + \sqrt{\frac{h(1-h)}{N}} \dot{W},
\end{equation}
where $N$ is the  number of particles involved and $\dot{W}$ is a normalized Gaussian white noise.
Contrary to Equation (\ref{equation F-KPP}), this equation with noise selects a single velocity $v_N$ for the front propagation.
A first approximation for Equation (\ref{equation F-KPP with noise}) consists in replacing the noise term by a deterministic cut-off, leading to the equation \cite{brunetderrida97,brunetderrida99}
\begin{equation} \label{equation F-KPP with cut-off}
\partial_t h 
= \frac{1}{2} \partial^2_x h + h(1-h) \1_{h \geq 1/N},
\end{equation}
which selects also a single velocity $v^\text{cutoff}_N$.
It has been conjectured in \cite{brunetderrida97,brunetderrida99} for the velocity $v^\text{cutoff}_N$ and in \cite{brunetderrida2001} for the velocity $v_N$ that, as $N$ tends to infinity,
\[
v_\text{min} - v_N 
\sim
v_\text{min} - v^\text{cutoff}_N 
\sim
\frac{\pi^2}{\sqrt{2} (\log N)^2},
\]
which is an extremely slow convergence.
These results have been given a rigorous mathematical proof in \cite{dpk2007} 
for the F-KPP equation with cut-off and more recently in \cite{mmq2011}
for the F-KPP equation with white noise.
Brunet and Derrida in \cite{brunetderrida97,brunetderrida99,brunetderrida2001} support their conjecture by studying directly a particular microscopic stochastic processes involving $N$ particles.
In the same way, Brunet, Derrida, Mueller and Munier \cite{bdmm2006-1} introduce three different branching-selection particle systems, in order to describe more precisely the velocity and the diffusion constant of the front for Equation (\ref{equation F-KPP with noise}).
Among these processes, they consider the $L$-BBM and the $N$-BBM. 
The latter is defined with a different selection rule: particles move according to Brownian motion, split after an exponential lifetime and, whenever the population size exceeds $N$, the lowest particle is killed.
They conjecture, among other things, that the asymptotic velocity of the $N$-BBM is $v_N$ and satisfies the more precise asymptotic expansion as $N\to\infty$:
\[
v_N 
=
\sqrt{2}
- 
\frac{\pi^2}{\sqrt{2} (\log N)^2}
\left(
1 -
\frac{6 \log \log N}{(\log N)^3} (1+\petito{1})
\right).
\]
This conjecture has not yet been mathematically proved, but some rigorous results have been obtained concerning the $N$-BBM or the $N$-branching random walk (its discrete time analog), 
in particular by Bérard and Gouéré \cite{berardgouere2010}, who showed the asymptotic behavior in $\pi^2/ \sqrt{2} (\log N)^2 (1 + \petito{1})$ for $\sqrt{2} - v_N$, 
but also by Durrett and Remenik \cite{durrettremenik2011}, Bérard and Maillard \cite{berardmaillard2014}, Mallein \cite{mallein2015arxiv} and Maillard \cite{maillard2015arxiv}.
For its part, the $L$-BBM has not yet been studied in the mathematical literature, but Brunet, Derrida, Mueller and Munier \cite{bdmm2006-1} conjectured that it behaves as the $N$-BBM by taking
\begin{equation} \label{equation L and N}
L
=
\frac{\log N}{\sqrt{2}},
\end{equation}
which means that, with (\ref{equation L and N}), the population size of the $L$-BBM is around $N$ and the $N$-BBM has approximately a width $L$.
It follows that the $L$-BBM must have an asymptotic velocity $v_L$ that satisfies, as $L \to \infty$,
\begin{equation} \label{equation dvlpt v_L}
\sqrt{2} - v_L \sim \frac{\pi^2}{2\sqrt{2}L^2},
\end{equation}
which is the result proved in this paper.
Moreover, some recent results of Berestycki, Berestycki and Schweinsberg \cite{bbs2013,bbs2011}, concerning BBM with absorption on a linear barrier and BBM in a strip,
suggested also that the asymptotic behavior (\ref{equation dvlpt v_L}) must hold.
Indeed, they show that for a BBM in the strip $(0,K)$ with drift $-\mu$ (it means that particles move according to Brownian motions with drift $-\mu$ and are killed by hitting 0 or $K$), the size of the population stays of the same order on a time scale of $K^3$ when
\[ 
\mu 
\coloneqq 
\sqrt{2 - \frac{\pi^2}{K^2}},
\]
see Proposition \ref{prop Z martingale} further.
Thus, if the fluctuations of the $L$-BBM around the deterministic speed $v_L$ are not too large (less than $\varepsilon L$ on a time scale of $L^3$), one can expect that 
\begin{equation} \label{oo}
\sqrt{2 - \frac{\pi^2}{(L - \varepsilon L)^2}} - \frac{\varepsilon L}{L^3}
\leq
v_L
\leq
\sqrt{2 - \frac{\pi^2}{(L + \varepsilon L)^2}} + \frac{\varepsilon L}{L^3},
\end{equation} 
for $\varepsilon >0$ and $L$ large enough, with comparison with BBM in strips $(0, L - \varepsilon L)$ and $(0, L + \varepsilon L)$,
and (\ref{equation dvlpt v_L}) follows from (\ref{oo}) by letting $\varepsilon \to 0$.
Actually, fluctuations of the $L$-BBM are believed to be of order $\log L$ on a time scale of $L^3$ and a precise understanding of them will probably lead to the next order in the asymptotic behavior of $v_L$.

\subsection{Proof overview and organization of the paper}

\begin{figure}[ht]
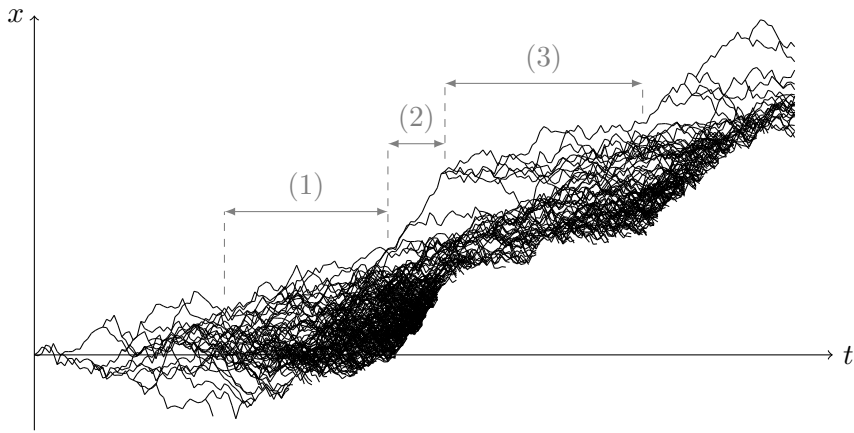
 
\centering 
 
\caption{Simulation between times 0 and 10 of an $L$-BBM with $L=5$ starting with a single particle at 0.
During the period (1), $\max X^L$ grows slowly so the number of particles becomes very large. 
Among all these particles, there is one that goes up very fast during the period (2) and it reduces drastically the population.
Therefore, after that, during period (3), $\max X^L$ grows once again slowly, involving an increase of the population size.} 
\label{figure L-BBM}
\end{figure}

One of the major difficulty in working with the $L$-BBM or the $N$-BBM is that they do not satisfy the branching property and, therefore, any form of many-to-one lemma (see Lemma \ref{lemma many-to-one}): the offspring of a particle at a time $t$ depends on the offspring of other particles alive at time $t$.
For this reason, Bérard and Gouéré \cite{berardgouere2010} compare the $N$-BRW with a branching random walk with absorption on a linear barrier in order to apply the precise results of Gantert, Hu and Shi \cite{ghs2011}.
In the same way, we come down here to results of Berestycki, Berestycki and Schweinsberg \cite{bbs2013} concerning branching Brownian motion in a strip (see Subsection \ref{subsection bbs}).
In both cases, the study is reduced to another process that satisfies the branching property and, thus, on which the work is easier.
However, the arguments used here are quite different from those of Bérard and Gouéré \cite{berardgouere2010}.
The main reason is the absence of a monotonicity property for the $L$-BBM, like the one used by Bérard and Gouéré (see Lemma 1 of \cite{berardgouere2010}): 
it does not seem to exist a coupling such that, when one of the initial particles of the $L$-BBM is removed, its maximum becomes stochastically smaller.
The proof of Proposition \ref{proposition existence} is based on the study of the return times to 1 for the population size.
These times delimit i.\@i.\@d.\@ pieces of the $L$-BBM, so, showing that they are sub-exponential, we can use the law of large numbers to prove the convergence of $\max X^L(t) / t$.
Although return times to 1 for the population size are sub-exponential, they are too large and, therefore, irrelevant for a more precise result.
Thus, to prove Theorem \ref{theorem}, we work instead with stopping times $(\tau_i)_{i\in\N}$ such that $\tau_{i+1} - \tau_i$ is shorter than $L^3$ and use on such a time interval a comparison with the BBM in a strip.
For the lower bound, we come down to a BBM in a strip that starts at time $\tau_i$ with exactly the same particles than the $L$-BBM and is then included in the $L$-BBM until time $\tau_{i+1}$.
Therefore, it is sufficient to show that this BBM in a strip goes up high enough between times $\tau_i$ and $\tau_{i+1}$.
It is done by using the monotonicity property of the BBM in a strip to consider the worst case with only a single initial particle at time $\tau_i$ and, then, applying  results of Berestycki, Berestycki and Schweinsberg \cite{bbs2013}.
In the same way, for the upper bound, we come down between times $\tau_i$ and $\tau_{i+1}$ to a BBM in a strip with more particles than the $L$-BBM and we show that it cannot rise too fast.
But here, we need to control the population size at time $\tau_i$: the bad cases happen when there are too many particles.
To this end, we use the following fact: a large population involves a quick increase of $\max X^L$ but, when the maximum of the $L$-BBM rises fast, many particle are killed by selection (see Figure \ref{figure L-BBM}).
This leads to the conclusion that $\max X^L$ cannot grow quickly during a too long period.
The paper is organized as follows. 
Section \ref{section results} introduces some useful results concerning the standard BBM and the BBM in a strip.
Then, Section \ref{section v_L} contains the proof of Proposition \ref{proposition existence}. 
Finally, the lower bound in Theorem \ref{theorem} is proved in Section \ref{section lower bound} and the upper bound in Section \ref{section upper bound}.
Throughout the paper, $C$ denotes a positive constant that does not depend on the parameters and can change from line to line.
For $x$ and $y$ real numbers, we set $x \wedge y \coloneqq \min(x,y)$ and $x \vee y \coloneqq \max(x,y)$.
For $f \colon \R \to \R$ and $g \colon \R \to \R_+^*$, we say that $f(x) \sim g(x)$ as $x\to\infty$ if $\lim_{x\to\infty} f(x)/g(x) = 1$
and that $f(x) = \petito{g(x)}$ as $x\to\infty$  if $\lim_{x\to\infty} f(x)/g(x) = 0$.
Lastly, we set $\N \coloneqq \{0,1,2,\dots \}$ and denote by $\sC([0,t],\R)$ the set of continuous functions from $[0,t]$ to $\R$.

\section{Some useful results} \label{section results}

\subsection{Standard branching Brownian motion} \label{subsection standard}

In this section, we present some classical results concerning branching Brownian motion without selection.
We assume here that there is initially one particle at 0, therefore we work under the probability $\P_{(0)}$ and the associated expectation $\E_{(0)}$, where $(0)$ denotes the configuration with a single particle at 0.
First of all, $(M(t))_{t\geq 0}$ under $\P_{(0)}$ is a Yule process (or pure birth process) with parameter 1, which means that it is a Markov process on $\N$ with transitions $i \to i+1$ at rate $i$ for all $i\in \N$.
Thus, $M(t)$ follows under $\P_{(0)}$ a geometric distribution with parameter $\e^{-t}$: for each $k \geq 1$, $\P_{(0)}(M(t)=k) = (1-\e^{-t})^{k-1} \e^{-t}$.
In particular, we have $\E_{(0)}[M(t)] = \e^t$ and the next result, often called in the literature the many-to-one lemma, follows.
\begin{lem}[Many-to-one lemma] \label{lemma many-to-one}
Let  $F \colon \sC([0,t],\R) \longrightarrow \R_+$ be a measurable function  and $(B_s)_{s \geq 0}$ denote a Brownian motion starting at 0. 
Then, we have
\[
\Eci{(0)}{\sum_{k=1}^{M(t)} F\left( (X_{k,t}(s))_{s\in [0,t]} \right)}
= \e^t \Ec{F\left( (B_s)_{s\in [0,t]} \right)},
\]
where we denote by $X_{k,t}(s)$ the position of the unique ancestor at time $s$ of $X_k(t)$.
\end{lem}
This lemma will be useful to compute the expectation of some functionals of the branching Brownian motion without selection. 
There is no similar result for the $L$-BBM because the number of particles of the $L$-BBM living at time $t$ is not independent of their trajectories on $[0,t]$.
We saw before that the first order for the asymptotic behavior of the position of the highest particle of the BBM is $\sqrt{2} t$, but we will need more accurate information on the extremal particle of the BBM.
The extremal particle has been a main topic in the study of the BBM, with in particular the seminal works of Bramson \cite{bramson78,bramson83} who shows the existence of a real random variable $W$ such that
\begin{equation} \label{cv en loi max BBM}
\max X(t) -m(t) 
\overset{\text{law}}{\underset{t\to\infty}{\longrightarrow}} W,
\end{equation}
where
\begin{equation} \label{definition m(t)}
m(t) 
\coloneqq
\sqrt{2} t - \frac{3}{2\sqrt{2}} \log t,
\end{equation}
and of Lalley and Sellke \cite{lalleysellke87} who describe the limit $W$ as a random mixture of Gumbel distributions.
We will need some further information on the trajectory leading to the extremal particle at time $t$.
The following proposition is an immediate consequence of convergence (\ref{cv en loi max BBM}) and of Theorem 2.5 of Arguin, Bovier and Kistler \cite{abk2011}.
It shows that there exists a high particle at time $t$ (above $m(t) - d$) whose trajectory does not go too far under the line $s \mapsto \frac{s}{t} m(t)$. 
Actually, Arguin, Bovier and Kistler \cite{abk2011} show that every high particle at time $t$ is likely to satisfy this property.
\begin{prop} \label{proposition ABK}
For all $\gamma >0$ and $\delta >0$, there exist $d>0$, $r>0$ and $t_0 >0$ large enough such that for all $t\geq t_0$,
\begin{align*}
\begin{split}
& \P_{(0)} \Bigl(
\exists k \in \llbracket 1,M(t) \rrbracket : 
X_k(t) \geq m(t) - d \\
& \hphantom{\P_{(0)} \Bigl(} \text{and } 
\forall s \in [0,t], 
X_{k,t}(s) \geq \frac{s}{t} m(t) 
	- r \vee \left(s^{\frac{1}{2}+\gamma} \wedge (t-s)^{\frac{1}{2}+\gamma} \right)
	 \Bigr)
\geq 
1 - \delta,
\end{split}
\end{align*}
where $X_{k,t}(s)$ denotes the position of the unique ancestor at time $s$ of $X_k(t)$.
\end{prop}

\subsection{Branching Brownian motion in a strip} \label{subsection bbs}

The branching Brownian motion with absorption on a linear barrier has been introduced by Kesten \cite{kesten1978}.
In this process, each particle moves according to a Brownian motion with drift $-\mu$, splits into two new particles after an exponentially distributed time and is killed when it reaches a non-positive position.
Kesten \cite{kesten1978} showed that the process survives with positive probability if and only if $\mu < \sqrt{2}$.
More recently, the branching Brownian motion with absorption has been studied by Harris, Harris and Kyprianou \cite{hhk2006}, by Harris and Harris \cite{harrisharris2007}
and, in the near-critical case, that is when $\mu \to \sqrt{2}$ while keeping $\mu < \sqrt{2}$, by Berestycki, Berestycki and Schweinsberg \cite{bbs2013,bbs2011}.
For recent results concerning the branching Brownian motion in a strip, see Harris, Hesse and Kyprianou \cite{hhk2016}.
We focus here on results of Berestycki, Berestycki and Schweinsberg \cite{bbs2013} concerning the branching Brownian motion in a strip, where particles move according to a Brownian motion with drift $-\mu$, split into two new particles after an exponentially distributed time and are killed when they come out of a fixed interval, that we will specify hereafter.
We fix a positive real $K$, we set 
\begin{equation} \label{equation mu K}
\mu \coloneqq \sqrt{2 - \frac{\pi^2}{K^2}}
\end{equation}
and we choose the interval $(0,K_A)$, where we set
\[ K_A \coloneqq K - \frac{A}{\sqrt{2}},\]
in order to keep notation similar to those of Berestycki, Berestycki and Schweinsberg \cite{bbs2013} (but with $K$ and $K_A$ instead of $L$ and $L_A$). 
We denote by $(\widetilde{X}^{K,K_A}_k(t), 1 \leq k \leq \widetilde{M}^{K,K_A}(t))_{t\geq 0}$ the BBM in the strip $(0,K_A)$ with drift $-\mu$, where $\mu$ is given by (\ref{equation mu K}).
When $A=0$ (that is $K = K_A$), we write $\widetilde{X}^K$ instead of $\widetilde{X}^{K,K_A}$, in order to simplify notation.
Berestycki, Berestycki and Schweinsberg \cite{bbs2013} introduce the following functional of the BBM in a strip, defined by
\begin{equation} \label{equation definition Z}
\widetilde{Z}^{K,K_A}(t) 
\coloneqq 
\sum_{k=1}^{\widetilde{M}^{K,K_A}(t)} \e^{\mu \widetilde{X}^{K,K_A}_k(t)}
\sin \left( \frac{\pi \widetilde{X}^{K,K_A}_k(t)}{K_A} \right),
\end{equation}
for each $t\geq 0$.
The random variable $\widetilde{Z}^{K,K_A}(t)$ estimates the size of the process and the process $(\widetilde{Z}^{K,K_A}(t))_{t\geq 0}$ has the advantage of being a martingale in the particular case where $K_A = K$, as stated in the following proposition (see Lemma 7 of \cite{bbs2013}). 
The $\sin$ function comes into play when one computes the density of a Brownian motion without drift started at $x\in (0,K)$ and killed when it hits 0 or $K$, by solving the heat equation on $(0,K)$ with zero as boundary conditions.
It brings also an exponential part that is then modified when one adds the drift and the branching events.
\begin{prop} \label{prop Z martingale}
The process $\e^{(1 - \mu^2/2 - \pi^2/2K_A^2)t} (\widetilde{Z}^{K,K_A}(t))_{t\geq 0}$ is a martingale.
In particular, when $A=0$, $(\widetilde{Z}^K(t))_{t\geq 0}$ is a martingale.
\end{prop}
So the drift $-\mu$ given by (\ref{equation mu K}) is exactly the right choice such that the size of the BBM in the strip $(0,K)$ does not degenerate.
When $A > 0$ (that is $K_A < K$), the population decreases exponentially and when $A < 0$ (that is $K_A > K$), the population explodes exponentially.
The following results concern the number of particles that are killed by hitting $K_A$ between times $0$ and $\theta K^3$.
This corresponds to Subsection 3.2 of Berestycki, Berestycki and Schweinsberg \cite{bbs2013} but with some modifications: here $\theta$ is not assumed small (rather, $\theta$ will be close to 1), we consider only the particular case $A=0$ and we do not make the same assumptions about the initial configuration.
This first proposition is a rewrite of Proposition 16 of Berestycki, Berestycki and Schweinsberg \cite{bbs2013} in the case $A=0$ and estimates the mean number of particles that are killed by hitting $K_A$ between times $0$ and $\theta K^3$ for any initial configuration.
\begin{prop} \label{prop estimation esperance de R}
Assume $A = 0$.
For some fixed $\theta > 0$, let $R$ be the number of particles that hit $K$ between times $0$ and $\theta K^3$ and $R'$ the number of particles that hit $K$ between times $K^{5/2}$ and $\theta K^3$.
Then, for all $\xi \in \cC$, $\P_\xi$-a.s.\@ we have the following inequalities as $K \to \infty$:
\begin{align*}
\Eci{\xi}{R}
& \leq
C \frac{\e^{-\mu K}}{K} \widetilde{V}^K(0)
+
2 \sqrt{2} \pi \theta K \e^{-\mu K}
\widetilde{Z}^K(0) (1 + \petito{1}), \\
\Eci{\xi}{R'}
& \geq
2 \sqrt{2} \pi \theta K \e^{-\mu K}
\widetilde{Z}^K(0) (1 + \petito{1}),
\end{align*}
where $C$ is a positive constant and $\widetilde{V}^K(0) \coloneqq \sum_{i=1}^{\widetilde{M}^K(0)} \widetilde{X}^K_i(0) \e^{\mu \widetilde{X}^K_i(0)}$.
\end{prop}
\begin{proof}
This result follows directly from the proof of Proposition 16 of Berestycki, Berestycki and Schweinsberg \cite{bbs2013}, by noting that $K = (\log N + 3 \log \log N)/\sqrt{2}$ and that the only place where they need $\theta$ to be small enough is their inequality (69), whose left-hand side is simply zero in the case $A=0$.
\end{proof}
The second proposition gives an upper bound for the second moment of $R$ and follows from the proof of Proposition 18 of \cite{bbs2013}.
We still work under the assumption $A=0$ but moreover with a single initial particle, therefore the initial configuration is denoted by $(x)$, where $x \in (0,K)$.
Here $x$ can be close to $K$, so the initial configuration does not necessarily satisfy the assumption of Proposition 18 of \cite{bbs2013}.
\begin{prop} \label{prop majoration esperance de R^2}
Assume $A = 0$
and suppose there is a single particle $x$ at time 0, where $x \in (0,K)$.
Let $R$ be the number of particles that hit $K$ between times $0$ and $K^3$.
Then, as $K\to\infty$,
\begin{align*}
\Eci{(x)}{R^2}
\leq
C \e^{-\mu K}	\e^{\mu x} 
\left( 1 + K \sin \left( \frac{\pi x}{K} \right) \right)
(1 + \petito{1}),
\end{align*}
where $C$ is a positive constant.
\end{prop}
\begin{proof}
In the same way as Berestycki, Berestycki and Schweinsberg \cite{bbs2013} in the proof of their Proposition 18, we write $R^2 = R + Y$, where $Y$ is the number of distinct pairs of particles that hit $K$ between times $0$ and $K^3$.
Here, as there is only a single initial particle, we have $Y = Y_2 = Y_2^x$ with notation of \cite{bbs2013}, which means that all particles that hit $K$ have the same ancestor at time 0.
From their bounds (76) and (77), it follows that
\[
\Eci{(x)}{Y}
\leq
C \e^{-\mu K} \e^{\mu x} (1 + \petito{1})
+
C \e^{-\mu K} K \e^{\mu x} 
	\sin \left( \frac{\pi x}{K} \right) (1 + \petito{1}),
\]
because they only use the assumption that $\theta \leq 1$ and we chose here $\theta = 1$.
Moreover, using Proposition \ref{prop estimation esperance de R}, we get
\begin{align*}
\Eci{(x)}{R}
& \leq
C \frac{\e^{-\mu K}}{K} x \e^{\mu x}
+
2 \sqrt{2} \pi K \e^{-\mu K}
\e^{\mu x} \sin \left( \frac{\pi x}{K} \right)
(1 + \petito{1}) \\
& \leq
C \e^{-\mu K}	\e^{\mu x} 
\left( 1 + K \sin \left( \frac{\pi x}{K} \right) \right)
(1 + \petito{1})
\end{align*}
and the result follows.
\end{proof}

\section{Existence of the asymptotic velocity}  \label{section v_L}

In this section, we prove Proposition \ref{proposition existence}.
We recall that the following proof is similar to the one presented by Derrida and Shi \cite{derridashi2016arxiv}.
However, the methods carried out here are very rough and cannot be used to get more elaborate results as Theorem \ref{theorem}, but we want to be ensured that the velocity $v_L$ of the $L$-branching Brownian motion is well defined. 
We work here with a fixed $L > 0$.
The strategy will be to show that the return time to 1 for the population size of the $L$-BBM is sub-exponential and then to apply the law of large numbers to the renewal structure obtained from the sequence of successive return times to 1 for the population size.
Foremost, we need the following coupling result that will be useful to show that there cannot be too many particles in the $L$-BBM.
\begin{lem} \label{lemma couplage}
Let $M$ be a positive integer and $\xi = (\xi_1,\dots,\xi_N)$ be a configuration, with $N \geq M$ and $\xi_1 \geq \dots \geq \xi_N$.
Then, there exist Brownian motions $B^1,\dots,B^M$ starting at $\xi_1, \dots, \xi_M$ respectively, that are mutually independent but not independent of the $L$-BBM, such that, for all $t \geq 0$, on the event $\{ \forall s \in [0,t], M^L(s) \geq M \}$, we have
\[
\max X^L(t) \geq \max_{1 \leq i \leq M} B^i_t
\]
$\P_\xi$-almost surely.
\end{lem}
\begin{proof}
We choose for $B^i$ the trajectory of the particle of the $L$-BBM that is at $\xi_i$ at time 0.
When this particle splits, we choose uniformly one of its children and $B^i$ continues by following the trajectory of this child.
We proceed like this until a particle followed by one of the Brownian motions $B^1,\dots,B^M$, say $B^i$, is killed by selection in the $L$-BBM.
Then, we distinguish two cases:
\begin{enumerate}
\item If there are at least $M$ particles in the $L$-BBM at this time denoted by $t$, then there is at least 1 particle that is not followed by one of the Brownian motions $B^1,\dots,B^M$, so we choose uniformly one of those particles and $B^i$ continues by following, from the position $B^i_t$, the trajectory of the chosen particle. 
We proceed then as before.
\item Otherwise, we lay the $L$-BBM aside and work only on the BBM without selection: each Brownian motion continues following its particle until it splits and then follows one of its child, and so on.
\end{enumerate}
Thus, $B^1,\dots,B^M$ are defined for all time and are independent Brownian motions.
On the event $\{ \forall s \in [0,t], M^L(s) \geq M \}$, we do not meet the case (ii) until time $t$ and so we have for all $s \in [0,t]$
\[
\max X^L(s) \geq \max_{1 \leq i \leq M} B^i_s,
\]
because $B^i$ at time $s$ is always lower than the particle it follows: when, in the case (i), $B^i$ changes the particle it follows, the new particle is necessarily above $B^i$, because the new particle is in the $L$-BBM and $B^i$ is at the low extremity of the $L$-BBM.
\end{proof}
We can now prove that the return time to 1 for the population size of the $L$-BBM is sub-exponential. 
For this purpose, we show with Lemma \ref{lemma couplage} that, if there is a very large number of particles in the population, then the maximum of the $L$-BBM can with high probability increase by $2L$ in a short time and, thus, a large proportion of particles is killed by selection.
\begin{prop} \label{prop temps de retour en 1}
Let $T \coloneqq \inf \{ t\geq 1 : M^L(t) = 1 \}$.
It exists a positive constant $c$ that depends on $L$ such that for all $\xi \in \cC$, 
$\Ppi{\xi}{T \leq 2} \geq c.$
\end{prop}
\begin{proof}
The strategy will be to show first the same result for $S_{M_0} \coloneqq \inf \{ t\geq 0 : M^L(t) \leq M_0 \}$, which is the return time under a fixed size $M_0$ for the population of the $L$-BBM, where $M_0$ will be chosen very large.
Let $M$ be a positive integer and $\xi = (\xi_1,\dots,\xi_N)$ a configuration, with $N \geq M$ and $\xi_1 \geq \dots \geq \xi_N$.
We can suppose with no loss of generality that $\xi_1 = 0$, because the law of $(M^L(t))_{t\geq0}$ is invariant under shift of the initial configuration, and that $\xi_N > -L$, otherwise the low particles will die instantly.
Let $a$ be a positive real number and $t \coloneqq aL^2 / \log N$. 
We have
\begin{align}
\begin{split}
\Ppi{\xi}{\forall s \in [0,t], M^L(s) \geq M}
&\leq
\Ppi{\xi}
{\forall s \in [0,t], M^L(s) \geq M, \max X^L(t) < 2L} \\
& \quad +
\Ppi{\xi}{M^L(t) \geq M, \max X^L(t) \geq 2L}.
\end{split} \label{ca}
\end{align}
The second term on the right-hand side of (\ref{ca}) can be bounded by noting that, on the event $\{ M^L(t) \geq M, \max X^L(t) \geq 2L \}$, there are at least $M$ particles of the BBM without selection that are above $L$.
Thus,
\begin{align}
\Ppi{\xi}{M^L(t) \geq M, \max X^L(t) \geq 2L}
& \leq
\Ppi{\xi}{\sum_{k=1}^{M(t)} \1_{X_k(t) \geq L} \geq M} \nonumber \\
& \leq \frac{1}{M} \Eci{\xi}{\sum_{k=1}^{M(t)} \1_{X_k(t) \geq L}} \nonumber \\
& = \frac{1}{M} \sum_{i=1}^N \e^t \Pp{B_t + \xi_i \geq L}, \label{cb}
\end{align}
where we used successively the Markov inequality and the many-to-one lemma (Lemma \ref{lemma many-to-one}) applied to the BBM starting at each $\xi_i$, with $(B_t)_{t\geq 0}$ denoting a Brownian motion starting at 0. 
With our choice for $t$ and recalling that $\xi_i \leq 0$, (\ref{cb}) is bounded from above by
\begin{align*}
\frac{N}{M} \e^{aL^2 / \log N} \Pp{ \sqrt{t} \cN(0,1) \geq L}
& \leq
\e^{aL^2 / \log N} \frac{N}{M} \frac{e^{- \log N / 2a}}{2}
=
\e^{aL^2 / \log N} \frac{N^{1-1/2a}}{2M},
\end{align*}
where $\cN(0,1)$ denotes the standard normal distribution and we use that 
$\P (\cN(0,1) \geq x) \leq \e^{-x^2/2}/2$ for all $x \geq 0$.
From now on, we set $M \coloneqq \lfloor N^{\lambda + 1-1/2a} \rfloor$, where $0 < \lambda < 1/2a$.
Thus, since $\e^{aL^2 / \log N} \to 1$ as $N\to\infty$, we proved that the second term on the right-hand side of (\ref{ca}) is smaller than $N^{-\lambda}$ for $N$ large enough.
We now deal with the first term on the right-hand side of (\ref{ca}). 
It can be bounded using Lemma \ref{lemma couplage}: we have, with $B^1,\dots,B^M$ independent Brownian motions starting at $\xi_1, \dots, \xi_M$ under $\P_\xi$,
\begin{align*}
\Ppi{\xi}
{\forall s \in [0,t], M^L(s) \geq M, \max X^L(t) < 2L}
&\leq 
\Ppi{\xi}{\max_{1 \leq i\leq M}  B^i(t) < 2L} \\
& \leq
\Pp{\sqrt{t}  \cN(0,1) -L <  2L}^M,
\end{align*}
since $\xi_i > -L$ for all $i$.
Using that $\Pp{\cN(0,1) \geq x} \geq \e^{-x^2/2}/2x$ for $x$ large enough and that $\ln(1-u) \leq -u$ for $u<1$, with our choice for $t$ and $M$, we get
\begin{align}
\Pp{\sqrt{t}  \cN(0,1) -L < 2L}^M
& =
\left( 
1- \Pp{\cN(0,1) \geq \frac{3 (\log N)^{1/2}}{a^{1/2}}} 
\right)^{\lfloor N^{\lambda + 1-1/2a} \rfloor} \nonumber \\
& \leq
\exp \left(
- \lfloor N^{\lambda + 1-1/2a} \rfloor 
\frac{N^{-9/2a} a^{1/2}}{6 (\log N)^{1/2}}
\right), \label{cc}
\end{align}
for $N$ large enough.
Thus, we want that $\lambda + 1- 10/2a > 0$, but we need $\lambda < 1/2a$ to have $M<N$.
So that a such $\lambda > 0$ exists, it is sufficient that $10/2a - 1 < 1/2a$, which means $a > 9/2$.
Then, we get $(\ref{cc}) \leq N^{-\lambda}$ for $N$ large enough.
Going back to (\ref{ca}), we showed that there exist $M_0$ large enough, $a>0$, $\lambda > 0$ and $0 < \mu < 1$ such that, for all $\xi = (\xi_1,\dots,\xi_N)$ with $N \geq M_0$,
\begin{equation} \label{cd}
\Ppi{\xi}{\exists s \in [0,aL^2/\log N] : M^L(s) \leq N^\mu} 
\geq 1 - N^{-\lambda}.
\end{equation}
We fix $N \geq M_0$ and we set 
$k \coloneqq \lceil \frac{\log \log N - \log \log M_0}{\log (1/ \mu)} \rceil$, 
so that $k$ is the single integer that satisfies 
$N^{\mu^k} \leq M_0 < N^{\mu^{k - 1}}$.
Then, by applying $k$ times the inequality (\ref{cd}) and the strong Markov property, we get
\begin{align}
\Ppi{\xi}{\exists s \in [0,a k L^2/\log N] : M^L(s) \leq N^{\mu^k}} 
&\geq 
(1 - N^{-\lambda}) (1 - (N^\mu)^{-\lambda}) \dotsm (1 - (N^{\mu^{k-1}})^{-\lambda}) 
\nonumber \\
& = 
\prod_{i=0}^{k-1} (1 - N^{-\lambda \mu^i})
\geq
\prod_{j=0}^{k-1} (1 - M_0^{-\lambda \mu^{-j}}), \label{ce}
\end{align}
by using that $N > M_0^{\mu^{-k+1}}$ and setting $j = k-i-1$.
It is easy to see that the product on the right-hand side of (\ref{ce}) converges to a positive limit as $k$ tends to infinity, if $M_0 >1$.
Moreover, we can choose $M_0$ large enough such that for all $N \geq M_0$, $a k L^2/\log N \leq 1$.
Thus, we have proved that there exist $c_1> 0$ and $M_0$ such that, for all $\xi = (\xi_1,\dots,\xi_N)$ with $N \geq M_0$, 
\begin{align}
\Ppi{\xi}{\exists s \in [0,1] : M^L(s) \leq M_0} \geq  c_1, \label{cf}
\end{align}
that is $\Ppi{\xi}{S_{M_0} \leq 1} \geq  c_1$ with $S_{M_0} \coloneqq \inf \{ t\geq 0 : M^L(t) \leq M_0 \}$.
Now, we consider an initial configuration $\xi = (\xi_1,\dots,\xi_N)$ but with $N \leq M_0$. 
We suppose $\xi_1 \geq \dots \geq \xi_N$.
If no particle splits on the time interval $[0,1]$, if the particle starting at $\xi_1$ stays above $\xi_1 - L/2$ on $[0,1]$ and reaches $\xi_1 + 2L$ at time 1 and if all other particles stay strictly under $\xi_1 + L/2$ on $[0,1]$, then at time 1 only the particle starting at $\xi_1$ is alive so $T=1$. 
Therefore, if $B$ denotes a Brownian motion starting at 0, we have
\begin{align*}
\Ppi{\xi}{T = 1} 
& \geq 
(e^{-1})^N
\Pp{\min_{s\in [0,1]} B_s \geq - \frac{L}{2}, B_t \geq 2L}
\Pp{\max_{s\in [0,1]} B_s < \frac{L}{2}}^{N-1}
\geq c_2,
\end{align*}
where $c_2 > 0$ is reached in the case $N = M_0$.
So we can conclude that $\forall \xi \in \cC$, $\Ppi{\xi}{T \leq 2} \geq c_1 c_2$, by using (\ref{cf}) and the strong Markov property at time $S_{M_0}$ in the case where there are more than $M_0$ particles in the initial configuration.
\end{proof}
The controls performed in the previous proof are very loose: with slightly more computation, one can see that $M_0$ needs to be larger than $\e^{b L^2}$ with some $b>0$, whereas most of the time $M^L(t)$ is of the order of $\e^{\sqrt{2} L}$.
Indeed, the return to 1 for the population size is a too infrequent event for a more accurate study of the $L$-BBM as in Sections \ref{section lower bound} and \ref{section upper bound}.
But it is sufficient here to prove the existence of the asymptotic velocity $v_L$ for the $L$-BBM.
\begin{proof}[Proof of Proposition \ref{proposition existence}]
Let $\xi \in \cC$ be a fixed configuration.
We set $T_0 \coloneqq 0$ and $T_{i+1} \coloneqq \inf \{ t\geq T_i + 1 : M^L(t) = 1 \}$ for each $i\in\N$.
Since for $i\geq 1$ there is a single particle in the $L$-BBM at time $T_i$, by the strong Markov property, $(T_{i+1} - T_i)_{i\geq1}$ and $(\max X^L(T_{i+1}) - \max X^L(T_i))_{i\geq1}$ are sequences of i.i.d.\@ random variables with the same laws as $T$ and $\max X^L(T)$ under $\P_{(0)}$ respectively.
Using Proposition \ref{prop temps de retour en 1} and the strong Markov property, we have for all $\xi\in\cC$ and $n \in \N$,
\[
\Ppi{\xi}{T \geq 2n} \leq (1-c)^n,
\]
which means that $T$ has a sub-exponential distribution under $\P_{\xi}$, so $T$ has finite moments and, in particular, $\P_\xi$-almost surely we have $T_1 < \infty$.
Therefore, by the law of large numbers, we get that $\P_\xi$-a.s.\@ $T_n/n \to \E_{(0)}[T]$ as $n \to \infty$.
It is clear that $\P_\xi$-a.s.\@ $\max X^L(T_1) < \infty$ but in order to apply the law of large numbers to the sequence $(\max X^L(T_{i+1}) - \max X^L(T_i))_{i\geq1}$, we need to check that $\E_{(0)}[\lvert \max X^L(T) \rvert]$ is finite. 
We are going to prove something stronger that will be useful afterwards in the proof, which is $\E_{(0)}[\zeta] < \infty$ where
\[
\zeta \coloneqq \max_{t\in [0,T]} \abs{\max X^L(t)}.
\]
For this, it is sufficient to prove that the function $a \mapsto \P_{(0)}(\zeta \geq a)$ is integrable on $\R_+$.
For all $a>0$, we have
\begin{align} \label{cg}
\Ppi{(0)}{\zeta \geq a}
& \leq 
\Ppi{(0)}{T > \sqrt{a}} 
+ 
\Ppi{(0)}{\exists t \in [0,\sqrt{a}] : 
\abs{\max X^L(t)} \geq a}.
\end{align}
Since $\E_{(0)}[T^2]$ is finite, $a \mapsto \P_{(0)}(T > \sqrt{a})$ is integrable on $\R_+$.
So we now have to deal with the second term on the right-hand side of (\ref{cg}).
By the coupling with the BBM without selection, it is bounded by
\begin{align} 
\Ppi{(0)}{\exists t \in [0,\sqrt{a}] : 
\abs{\max X(t)} \geq a} 
& \leq
2 \Ppi{(0)}{\exists t \in [0,\sqrt{a}] : 
\max X(t) \geq a} \nonumber \\
&\leq
2 \Eci{(0)}{\sum_{k=1}^{M(\sqrt{a})} 
\1_{\exists t \in [0,\sqrt{a}] : 
X_{k,\sqrt{a}}(t) \geq a}}, \label{ch}
\end{align}
where we denote by $X_{k,\sqrt{a}}(t)$ the position of the unique ancestor at time $t$ of $X_k(\sqrt{a})$.
Using the many-to-one lemma (Lemma \ref{lemma many-to-one}), (\ref{ch}) is equal to
\begin{align*} 
2 \e^{\sqrt{a}} 
\Pp{\exists t \in [0,\sqrt{a}] : B_t \geq a} 
& =
4 \e^{\sqrt{a}} \Pp{a^{1/4} \cN(0,1) \geq a} 
\leq
2 \e^{\sqrt{a}} \e^{- a^{3/2}/2},
\end{align*}
which is an integrable function of $a$.
This concludes the proof of the fact that $\E_{(0)}[\zeta]$ is finite.
In particular, we can apply the law of large numbers to the sequence $(\max X^L(T_{i+1}) - \max X^L(T_i))_{i\geq1}$ and get
that $\P_\xi$-a.s.\@ $\max X^L(T_n)/n \to \E_{(0)}[\max X^L(T)]$ as $n \to \infty$.
Thus, we have the convergence 
\begin{equation} \label{ci}
\frac{\max X^L(T_n)}{T_n}
\underset{n\to\infty}{\longrightarrow}
\frac{\E_{(0)}[\max X^L(T)]}{\E_{(0)}[T]}
\eqqcolon
v_L
\end{equation}
$\P_\xi$-almost surely.
For each $t \geq 0$, let $n_t$ be the integer such that $T_{n_t} \leq t < T_{n_t +1}$. 
It suffices now to show that $\P_\xi$-a.s.\@ $\max X^L(t)/t - \max X^L(T_{n_t})/T_{n_t} \to 0$\@ as $t\to\infty$.
We have
\begin{align}
\abs{\frac{\max X^L(t)}{t} - \frac{\max X^L(T_{n_t})}{T_{n_t}}}
& \leq
\abs{\frac{\max X^L(t) - \max X^L(T_{n_t})}{t}}
+
\abs{\frac{\max X^L(T_{n_t})}{T_{n_t}}} \left( \frac{t-T_{n_t}}{t} \right) \nonumber \\
& \leq
\frac{\zeta_{n_t}}{t}
+
\abs{\frac{\max X^L(T_{n_t})}{T_{n_t}}} \left( 1 - \frac{T_{n_t}}{t} \right), \label{cj}
\end{align}
where, for $n\in \N$, we set 
\[ 
\zeta_n 
\coloneqq 
\max_{t\in [T_n,T_{n+1}]} \abs{\max X^L(t) - \max X^L(T_n)}.
\]
Since $T_{n_t}/n_t \leq t/n_t < T_{n_t + 1}/n_t$, we get that $\P_\xi$-a.s.\@ $t/n_t \to \E_{(0)}[T]$ as $n \to \infty$, so $T_{n_t}/t \to 1$ and combining with (\ref{ci}), we deduce that $\P_\xi$-a.s.\@ the second term on the right-hand side of (\ref{cj}) tends to 0.
We now have to deal with the first term.
Since $(\zeta_n)_{n\geq 1}$ is a sequence of i.i.d.\@ random variables with the same law as $\zeta$ under $\P_{(0)}$ and $\E_{(0)}[\zeta]$ is finite, we have by the law of large numbers $\P_\xi$-a.s.\@ $\xi_n / n \to 0$ as $n\to\infty$
and so $\P_\xi$-a.s.\@ $\zeta_{n_t} / t \to 0$ as $t\to\infty$.
It concludes the proof of Proposition \ref{proposition existence}.
\end{proof}

\section{Lower bound for \texorpdfstring{$v_L$}{vL}} \label{section lower bound}

In this section, we fix $0< \varepsilon <1$ and consider all processes with drift $-\mu$, where
\[ \mu \coloneqq \sqrt{2 - \frac{\pi^2}{(1-\varepsilon)^2 L^2}},\]
which means that, for each $\xi \in \cC$, under $\P_\xi$,
$(X_k(t), 1 \leq k \leq M(t))_{t\geq 0}$ is a BBM without selection, with drift $-\mu$ and starting from the configuration $\xi$ and $(X^L_k(t), 1 \leq k \leq M^L(t))_{t\geq 0}$ is the associated $L$-BBM with drift $-\mu$.
However, $v_L$ still denotes the asymptotic velocity of $L$-BBM without drift, so Proposition \ref{proposition existence} shows that, for all $L>0$ and all $\xi \in \cC$,
\begin{equation} 
\frac{\max X^L(t)}{t}
\underset{t\to\infty}{\longrightarrow}
v_L - \mu
\end{equation}
$\P_\xi$-almost surely. 
Actually, the aim is to show that for $L$ large enough $\lim \max X^L(t) / t \geq 0$  $\P_{(0)}$-a.s.\@ so that we can conclude that $v_L \geq \mu$ and the lower bound follows by letting $\varepsilon \to 0$.

\subsection{Proof of the lower bound}

\begin{figure}[ht]
\centering
\begin{tikzpicture}
\fill[gray!30] (0,4) -- (7,4) -- (7,1) -- (0,1) -- cycle;

\node (d) at (0.5,0.36) {};
\node (f) at (-0.6/0.8,0.16) {};
\node[gray,text width=3cm,text centered] (t) at (-2,0.16) {Killing barrier of the $L$-BBM};
\draw[->,>=latex,gray] (d) to (f.east);
\draw[smooth,blue!70] plot coordinates 
{(0,3.5-3) (0.2,3.52-3) (0.5,3.2-3) (0.7,3.15-3) (1,2.8-3) (1.3,2.9-3) (1.5,2.9-3) (1.9,3.2-3) (2.05,2.925-3)};
\draw[smooth,blue!70] plot coordinates 
{(2.05,2.925-3) (2.2,3.1-3) (2.5,3.15-3) (2.8,3.5-3) (3.2,3.8-3) (3.5,3.5-3) (3.7,3.6-3) (3.9,3.6-3) (4.2,3.8-3) (4.5,4-3)};

\draw[->,>=latex] (-0.4,0) -- (8,0) node[right]{$t$};
\draw[->,>=latex] (0,-0.08) -- (0,5) node[left]{$x$};

\draw (0,-0.08) node[below]{$\tau_i$};
\draw (7,0.08) -- (7,-0.08) node[below]{$\tau_i + L^3$};
\draw (4.5,0.08) -- (4.5,-0.08) node[below]{$\tau_{i+1}$};

\draw (-0.08,3.5) node[left]{$\max X^L (\tau_i)$};
\draw (-0.08,4.05) node[left]{$\max X^L (\tau_i) + 1$};
\draw (-0.08,1) node[left]{$\max X^L (\tau_i) - (L-1)$};

\draw[thick] (-0.08,4) -- (7,4);
\draw[thick] (-0.08,1) -- (7,1);
\draw[thick] (7,4) -- (7,1);
\draw[<->,>=latex] (7.3,1) -- (7.3,4);
\draw (7.3,2.5) node[right]{$L$};

\draw[dotted] (4.5,0) -- (4.5,4); 
\draw[dotted] (7,0) -- (7,1);

\draw[smooth] plot coordinates 
{(0,3.5) (0.2,3.52) (0.5,3.2) (0.7,3.15) (1,2.8) (1.3,2.9) (1.5,2.9) (1.9,3.2) (2.05,2.925) (2.2,2.65) (2.4,2.55) (2.8,2) (3.2,2.05) (3.6,1.2) (3.9,1.1) (4.2,0.8)};
\draw[smooth] plot coordinates 
{(1.5,2.9) (2,2.75) (2.5,2.2) (3,2.7) (3.3,2.75) (3.9,3.3) (4.25,3.3) (4.5,3)};

\draw[smooth] plot coordinates 
{(0,2.7) (0.2,2.65) (0.5,2.4) (0.7,2.4) (1,2.2) (1.3,2.25) (1.6,2.7) (1.9,2.75) (2.05,2.925) (2.2,3.1) (2.5,3.15) (2.8,3.5) (3.2,3.8) (3.5,3.5) (3.7,3.6) (3.9,3.6) (4.2,3.8) (4.5,4)};
\draw[smooth] plot coordinates 
{(2.8,3.5) (3,3.45) (3.2,3.1) (3.4,3.1) (3.7,2.65) (3.9,2.7) (4.3,2.4) (4.5,2.4)};

\draw[smooth] plot coordinates 
{(0,1.9) (0.2,1.8) (0.5,2.1) (0.8,1.9) (1,1.9) (1.3,1.5) (1.5,1.5) (2,1.2) (2.3,1.4) (2.5,1.4) (2.8,1.6) (3.3,1.4) (3.7,1.45) (4,1.8) (4.3,1.75) (4.5,1.8)};
\draw[smooth] plot coordinates 
{(1,1.9) (1.3,1.9) (1.5,1.8) (1.8,1.45) (2,1.4) (2.3,0.9) (2.6,0.8) (2.8,0.5)};
\draw[smooth] plot coordinates 
{(3.7,1.45) (4,1.3) (4.3,1.35) (4.5,1.2)};

\draw[smooth] plot coordinates 
{(0,1.3) (0.3,1.4) (0.5,1.1) (0.7,1.05) (1,0.8) (1.3,0.9) (1.5,0.8) (1.7,0.4) (1.9,0.35) (2.3,0.4) (2.5,3.15-3)};
\draw[smooth] plot coordinates 
{(1.9,0.35) (2.2,3.1-3)};
\end{tikzpicture}
\caption{Representation of a $L$-BBM between times $\tau_i$ and $\tau_{i+1}$. 
By definition, $\tau_{i+1}$ is the time where $t \mapsto \max X^L (t)$ leaves the gray area. 
We are here on the event $A_i^c$ because $t \mapsto \max X^L (t)$ leaves the area from above.
Note that the killing barrier of $X^L$ (drawn in blue) stays below $\max X^L (\tau_i) - (L-1)$.} 
\label{figure definition tau_i lower bound}
\end{figure}
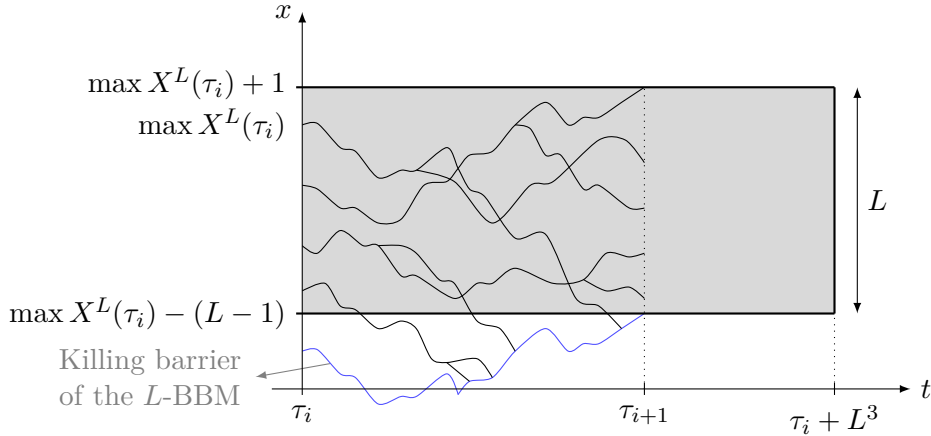

In this subsection, we prove the lower bound in Theorem \ref{theorem}, by postponing to the next subsection the proof of a proposition.
The strategy is to study the $L$-BBM on time intervals of length at most $L^3$ associated to a sequence of stopping times $(\tau_i)_{i\in\N}$ defined by $\tau_0 \coloneqq 0$ and for each $i \in \N$, 
\[ 
\tau_{i+1} 
\coloneqq 
(\tau_i + L^3) \wedge
\inf \enstq{t \geq \tau_i}
	{\max X^L(t) - \max X^L(\tau_i) \notin (-(L-1),1)}.
\]
We also define the event $A_i \coloneqq \{ \max X^L(\tau_{i+1}) < \max X^L(\tau_i) + 1 \}$ (see Figure \ref{figure definition tau_i lower bound} for an illustration of these definitions).
The event $A_i$ is a ``bad'' event, because on the event $A_i^c$ the position of the highest particle of the $L$-BBM goes up between times $\tau_i$ and $\tau_{i+1}$.
The following proposition shows that the event $A_i$ is very unlikely, regardless of the configuration of the $L$-BBM at time $\tau_i$. 
Its proof is postponed to Subsection \ref{subsection control P(A) lower bound}.
\begin{prop} \label{prop control P(A) lower bound}
Let $h(L) \coloneqq \sup_{\xi \in \cC} \P_\xi(A_0)$. 
Then, as $L\to\infty$, we have $h(L) = \petito{1/L}$.
\end{prop}
Using Proposition \ref{prop control P(A) lower bound}, we can now conclude the proof of the lower bound.
Let $K_n \coloneqq \sum_{i=0}^{n-1} \1_{A_i}$ be the number of ``bad'' events $A_i$ that happen before time $\tau_n$.
On the event $A_i^c$, we have $\max X^L(\tau_{i+1}) - \max X^L(\tau_i) = 1$ and, on the event $A_i$, we only have $\max X^L(\tau_{i+1}) - \max X^L(\tau_i) \geq - (L-1)$ and, therefore, we get
\[
\max X^L(\tau_n) - \max X^L(0)
\geq
(n-K_n) - (L-1) K_n
=
n - L K_n
\]
and, thus,
\begin{equation} \label{da}
\frac{\max X^L(\tau_n) - \max X^L(0)}{\tau_n}
\geq
\frac{n}{\tau_n} \left( 1 - L \frac{K_n}{n} \right).
\end{equation}
We now want to prove that the right-hand side of (\ref{da}) is non-negative as $n$ tends to infinity.
For this, we need to control $K_n$. 
Using Proposition \ref{prop control P(A) lower bound} and the strong Markov property, we get that, for each $0 \leq k \leq n$,
$
\P_{(0)}(K_n \geq k)
\leq 
\binom{n}{k} h(L)^k
$
and so
\[
\Eci{(0)}{\e^{K_n}}
\leq \sum_{k=0}^n \binom{n}{k} h(L)^k \e^k
= (1 + h(L) \e)^n
\leq \e^{n h(L) \e}.
\]
Then, applying the Markov inequality, we get
\[
\Ppi{(0)}{K_n \geq 3nh(L)} 
\leq \e^{-3nh(L)} \Eci{(0)}{\e^{K_n}}
\leq \e^{-nh(L) (3-\e)}
\]
which is summable, so the Borel-Cantelli lemma implies that $\P_{(0)}$-a.s.
$\limsup_{n\to \infty} \frac{K_n}{n} \leq 3h(L)$.
Moreover, using Proposition \ref{prop control P(A) lower bound}, we get that for $L$ large enough $1 - 3 h(L) L \geq 1/2$, so we conclude that for $L$ large enough,
\begin{equation} \label{db}
\liminf_{n\to \infty} \frac{n}{\tau_n} \left( 1 - L \frac{K_n}{n} \right) \geq 0.
\end{equation}
Assume for now that $\tau_n$ tends to infinity $\P_{(0)}$-almost surely as $n \to \infty$.
Then, the left-hand side of (\ref{da}) converges $\P_{(0)}$-almost surely to $v_L - \mu$ as $n \to \infty$, so with (\ref{db}) we get $v_L - \mu \geq 0$ for $L$ large enough, that is
\[
v_L \geq \sqrt{2} - \frac{\pi^2}{2\sqrt{2} (1-\varepsilon)^2 L^2} + \petito{\frac{1}{L^2}}
\]
and the lower bound in Theorem \ref{theorem} follows by letting $\varepsilon \to 0$.
It remains to show that $\tau_n$ tends to infinity $\P_{(0)}$-almost surely.
On the event $\{ (\tau_n)_{n \geq 0} \text{ is bounded} \}$, 
$\tau_n$ tends to a finite limit $\tau_\infty$ as $n \to \infty$, so it follows from (\ref{da}) that $\P_{(0)}$-almost surely
\[
\frac{\max X^L(\tau_\infty)}{\tau_\infty} = \infty,
\]
for $L$ large enough such that $1 - 3 h(L) L \geq 1/2$. 
But $\P_{(0)}$-almost surely, for all $t\geq 0$, $\max X^L(t)/t < \infty$, 
so $\P_{(0)} ((\tau_n)_{n \geq 0} \text{ is bounded}) = 0$.
It concludes the proof of the lower bound in Theorem \ref{theorem}.

\subsection{Proof of Proposition \ref{prop control P(A) lower bound}} 
\label{subsection control P(A) lower bound}

In this subsection, we prove Proposition \ref{prop control P(A) lower bound}, that is we show that the event $A_0^c$ is very likely. 
For this, we compare the $L$-BBM with a BBM in a strip that has less particles than the $L$-BBM between times 0 and $\tau_1$. 
Thus, if a particle of the BBM in a strip reaches $\max X^L(0) + 1$, then it is also the case for the $L$-BBM.
To prove that there is a particle of the BBM in a strip reaching $\max X^L(0) + 1$ with high probability, we will use Proposition \ref{prop estimation esperance de R}, Proposition \ref{prop majoration esperance de R^2} and the following lemma that allows to have many independent BBM in a strip trying to reach $\max X^L(0) + 1$ instead of only one.
\begin{lem} \label{lemma branchement du BBM}
Assume that $X$ is a BBM with drift $-\mu$ where $\mu < \sqrt{2}$. Then, for each $\alpha > 0$, for $L$ large enough, we have
\[
\Ppi{(0)}{M(L^\alpha) < L \text{ or } 
\exists t \in [0,L^\alpha]: \min X(t) \leq -4L^\alpha}
\leq
2L \e^{-L^\alpha}.
\]
\end{lem}
Note that this lemma concerns branching Brownian motion without selection and is stated for a large choice of drift $-\mu$ so that it can be applied here and also in Section \ref{section upper bound}.
\begin{proof}
Since $M(L^\alpha)$ follows a geometric distribution with parameter $\e^{-L^\alpha}$, we get
\begin{equation} \label{ea}
\Pp{M(L^\alpha) < L}
=
\sum_{k=1}^{L-1} \e^{-L^\alpha} \left( 1-\e^{-L^\alpha} \right)^{k-1}
=
1- \left( 1-\e^{-L^\alpha} \right)^{L-1}
\underset{L\to\infty}{\sim}
L \e^{-L^\alpha}.
\end{equation}
Moreover, applying the many-to-one lemma (Lemma \ref{lemma many-to-one}) without forgetting the drift $-\mu$, we have, with $(B_t)_{t\geq0}$ a Brownian motion,
\begin{align}
\Pp{\exists t \in [0,L^\alpha] : \min X(t) \leq -4L^\alpha}
& \leq 
\e^{L^\alpha} 
\Pp{\exists t \in [0,L^\alpha] : B_t - \mu t \leq -4L^\alpha} 
	\nonumber \\
& \leq 
\e^{L^\alpha} 
\Pp{\exists t \in [0,L^\alpha] : B_t \leq -(4 - \mu) L^\alpha} 
	\nonumber \\
& = 
\e^{L^\alpha} 2 \Pp{L^{\alpha/2} \cN(0,1) > (4 - \mu) L^\alpha} \nonumber \\
& \leq
\e^{L^\alpha - (4 - \mu)^2 L^\alpha /2}. \label{eb}
\end{align}
The result follows from (\ref{ea}) and (\ref{eb}) with $\mu \leq \sqrt{2}$.
\end{proof}
We now have enough tools to prove Proposition \ref{prop control P(A) lower bound}. 
Actually, the control will be much more accurate than needed.
\begin{proof}[Proof of Proposition \ref{prop control P(A) lower bound}]
Let $\xi \in \cC$ be a configuration.
Note first that $\P_\xi (A_0)$ is invariant under shift of $\xi$, so we can assume without loss of generality that $\max \xi = L-1$.
Then, $\P_\xi$-a.s.\@ $\tau_1 = L^3 \wedge 
\inf \{ t \geq 0: \max X^L(t) \notin (0,L) \}$
and $A_0$ is $\P_\xi$-a.s.\@ equal to the event ``no particle of the $L$-BBM reaches $L$ before reaching 0 on the time interval $[0,L^3]$''.
Moreover, $\P_\xi$-a.s.\@ on the time interval $[0,\tau_1]$ the killing barrier of the $L$-BBM stays under 0 (see Figure \ref{figure definition tau_i lower bound}), so we have the inclusion
\begin{equation} \label{ec}
\forall t \in [0,\tau_1], \widetilde{X}^{K,L}(t) \subset X^L(t),
\end{equation}
where we set $K \coloneqq (1-\varepsilon) L$ so that (\ref{equation mu K}) is satisfied and $A \coloneqq -\varepsilon L$ so that $K_A = L$: 
thus, under $\P_\xi$, $(\widetilde{X}^{K,L}_k(t), 1 \leq k \leq \widetilde{M}^{K,L}(t))_{t\geq 0}$ is a BBM in the strip $(0,L)$ with drift $-\mu$ starting from the configuration $\xi$.
Let $C_0$ denote the event ``no particle of $\widetilde{X}^{K,L}$ reaches $L$ on the time interval $[0,L^3]$'', then it follows from (\ref{ec}) that $\P_\xi$-a.s.\@ $A_0 \subset C_0$.
Moreover, the BBM in a strip satisfies the branching property : the offspring of a single particle at $x$ at time 0 is independent of the offspring of other initial particles and follows the law of a BBM in a strip under $\P_x$.
So we have $\P_\xi (C_0) \leq \P_{(L-1)}(C_0)$, by keeping only the offspring of the highest initial particle.
Thus, we get $h(L) \leq \P_{(L-1)}(C_0)$ because it does not depend any more on the initial configuration.
Now, our aim is to give an upper bound for $\P_{(L-1)}(C_0)$. 
For this purpose, we will first use Lemma \ref{lemma branchement du BBM} to have at least $L$ particles above $L-1 - 4L^\alpha$ after a short time $L^\alpha$ and then apply Propositions \ref{prop estimation esperance de R} and \ref{prop majoration esperance de R^2} to show that each particle at time $L^\alpha$ has a descendant that reaches $L$ before time $L^3$ with a positive probability that does not depend on $L$.
We fix $0 < \alpha < 1/2$.
Using Lemma \ref{lemma branchement du BBM}, we get, for $L$ large enough such that $L-1-4L^\alpha >0$,
\begin{align}
& \P_{(L-1)}(C_0) \nonumber \\
& \leq 
2L \e^{-L^\alpha}
+
\Ppi{(L-1)}
{ \{ M(L^\alpha) \geq L \} 
\cap \{ \forall t \in [0,L^\alpha], \min X(t) > L-1-4L^\alpha \}
\cap C_0} \nonumber \\
& \leq 
2L \e^{-L^\alpha}
+
\Ppi{(L-1)}
{ \{ \widetilde{M}^{K,L}(L^\alpha) \geq L \} 
\cap \{ \min \widetilde{X}^{K,L}(L^\alpha) > L-1-4L^\alpha \}
\cap C_0}, \label{ed}
\end{align}
because on the event $\{ \forall t \in [0,L^\alpha], \min X(t) > -4L^\alpha \} \cap C_0$, no particle of the BBM in a strip is killed between times 0 and $L^\alpha$, so $\widetilde{M}(L^\alpha) = M(L^\alpha)$.
Applying the branching property at time $L^\alpha$, we bound from above the second term of (\ref{ed}) by
\begin{equation} \label{ee}
\Eci{(L-1)}{ 
\1_{\widetilde{M}^{K,L}(L^\alpha) \geq L} 
\1_{\min \widetilde{X}^{K,L}(L^\alpha) > L-1-4L^\alpha}
\prod_{k=1}^{\widetilde{M}^{K,L}(L^\alpha)} \Ppi{(\widetilde{X}^{K,L}_k(L^\alpha))}{C_0'}}, 
\end{equation}
where $C_0'$ denotes the event ``no particle of $\widetilde{X}^{K,L}$ reaches $L$ on the time interval $[0,L^3 - L^\alpha]$''.
Note that the function $x \in (0,L) \mapsto \P_{(x)}(C_0')$ is nondecreasing\footnote{It follows from the fact that $\P_{(x)}(C_0') = \P_{(x)}( \text{no particle of } \hat{X} \text{ reaches } L \text{ on } [0,L^3 - L^\alpha])$, where $\hat{X}$ denotes the BBM with drift $-\mu$, with absorption at 0 and with a single initial particle at $x$ under $\P_{(x)}$.}.
So it follows from (\ref{ee}) that 
\begin{equation} \label{ef}
h(L) 
\leq 
2L \e^{-L^\alpha}
+
\Ppi{(L-1-4L^\alpha)}{C_0'}^L,
\end{equation}
for $L$ large enough.
Our aim is now to control $\P_{(L-1-4L^\alpha)}(C_0')$.
We set $K' \coloneqq L - 4L^\alpha$ and
\[
\mu' \coloneqq \sqrt{2 - \frac{\pi^2}{(L - 4L^\alpha)^2}},
\]
the drift associated to $K'$ (see equation (\ref{equation mu K})).
Moreover, we define a new process $\bar{X}$ from the standard BBM $X$ with drift $-\mu$ by killing particles that go below $t \mapsto (\mu' - \mu) t$ or above $t \mapsto K' + (\mu' - \mu) t$ (see Figure \ref{figure changement de drift}).
\begin{figure}[ht]
\centering
\begin{tikzpicture}
\draw[->,>=latex] (-0.08,0) -- (7,0) node[right]{$t$};
\draw[->,>=latex] (0,-0.08) -- (0,5.3) node[left]{$x$};
\draw (0,-0.08) node[below left]{$0$};

\draw (6,4) -- (-0.08,4) node[left]{$L$};
\draw (0.08,2.8) -- (-0.08,2.8) node[left]{$(1-\varepsilon) L = K$};

\draw[thick] (6,5) -- (-0.08,3.6-0.08*14/60) node[left]{$L - 4L^\alpha = K'$};
\draw[thick] (6,1.4) -- (0,0);

\draw[dotted] (12/7,4) -- (12/7,0);
\draw (12/7,0.08) -- (12/7,-0.08) node[below]{$\frac{L-K'}{\mu' - \mu}$};
\draw[dotted] (6,5) -- (6,0);
\draw (6,0.08) -- (6,-0.08) node[below]{$L^3 - L^\alpha$};

\draw[smooth] plot coordinates 
{(0,3.45) (0.2,3.45) (0.5,3.2) (1,3.3) (1.5,2.9) (1.9,3) (2.6,2) (3,2.05) (3.4,1.2) (3.7,1.1) (3.9,14*3.9/60)};

\draw[smooth,dashed] plot coordinates 
{(3.9+0.05,14*3.9/60-0.07) (4.1,0.7) (4.4,0.7) (4.8,0.4) (5.2,0.45) (5.5,0.6) (6,0.4)};
\draw[smooth,dashed] plot coordinates 
{(4.8,0.4) (5,0.3) (5.15,0)};
\draw[smooth] plot coordinates 
{(3,2.05) (3.4,2.1) (3.7,1.9) (4,1.95) (4.5,1.3) (4.7,1.3) (4.9,14*4.9/60)};
\draw[smooth,dashed] plot coordinates 
{(4.9+0.05,14*4.9/60-0.06) (5.1,1) (5.4,1.1) (5.7,1.4) (6,1.6)};

\draw[smooth] plot coordinates 
{(1,3.3) (1.5,3.5) (2,3.3) (2.5,3.25) (2.8,3.5) (3.2,3.8) (3.9,3.6) (4.2,3.8) (4.5,3.9) (5,3.7) (5.5,3.3) (6,3.2)};
\draw[smooth] plot coordinates 
{(2.8,3.5) (3.2,3.4) (3.9,3.1) (4.2,3.15) (4.5,2.7) (5,2.6) (5.5,2.65) (6,2.4)};
\draw[smooth] plot coordinates 
{(5,3.7) (5.3,3.8) (5.7,3.65) (6,3.6)};

\draw[smooth] plot coordinates 
{(0.2,3.45) (0.4,3.5) (0.6,3.6+14*0.6/60)};
\draw[smooth,dashed] plot coordinates 
{(0.6+0.05,3.6+14*0.6/60+0.08) (0.8,3.9) (1.2,3.6) (1.4,3.6) (1.8,3.2) (2.4,3) (2.6,3.05) (3,2.7) (3.4,2.9) (3.6,3) (3.9,2.9) (4.3,3) (4.6,3.05) (5,2.9) (5.3,2.85) (5.6,3.15) (6,3.4)};
\draw[smooth,dashed] plot coordinates 
{(3,2.7) (3.5,2.5) (3.8,2.55) (4.2,2.3) (4.7,2.2) (5.3,1.7) (6,1.9)};

\draw[->,>=latex,gray] (4,3.6 + 14/60*4) -- (3,5) node[left]{slope $\mu' - \mu$};

\end{tikzpicture}
\caption{Representation of the coupled systems $\bar{X}$ (full line) and $\widetilde{X}^{K,L}$ (dashed line) starting with a single initial particle at $K'-1$, on the event $C_0'$. 
The two thick straight lines of slope $\mu' - \mu$ are the killing barriers that define $\bar{X}$.} 
\label{figure changement de drift}
\end{figure}
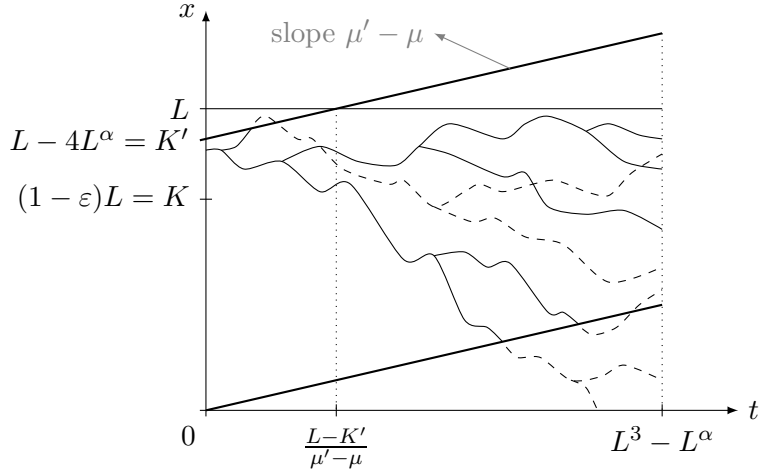
With this coupling, on the event $C_0'$, we have the inclusion $\forall t \in [0,L^3 - L^\alpha], \bar{X}(t) \subset \widetilde{X}^{K,L}(t)$ and it follows that
\begin{align*}
C_0'
& \subset
\left\{
\text{no particle of } \bar{X} \text{ reaches } L \text{ on time interval } 
\left[ \frac{L-K'}{\mu' - \mu}, L^3 - L^\alpha \right]
\right\} \\
& \subset
\left\{
\text{no particle of } \bar{X} \text{ reaches } 
t \mapsto K' + (\mu' - \mu) t 
\text{ on } 
\left[ \frac{L-K'}{\mu' - \mu}, L^3 - L^\alpha \right]
\right\},
\end{align*}
where the last event means that $\forall t \in [ (L-K')/(\mu' - \mu), L^3 - L^\alpha ], \max \bar{X}(t) < K' + (\mu' - \mu) t$.
Then, recalling that $\widetilde{X}^{K'}$ denotes the BBM in the strip $(0,K')$ with drift $-\mu'$, note that 
\begin{align*}
& \Ppi{(L-1-4L^\alpha)}{
\text{no particle of } \bar{X} \text{ reaches } 
t \mapsto K' + (\mu' - \mu) t 
\text{ on } 
\left[ \frac{L-K'}{\mu' - \mu}, L^3 - L^\alpha \right]
} \\
& \leq
\Ppi{(L-1-4L^\alpha)}{
\text{no particle of } \widetilde{X}^{K'} \text{ reaches } K' \text{ on } 
\left[ \frac{L-K'}{\mu' - \mu}, L^3 - L^\alpha \right]
} \\
& \leq
\Ppi{(L-1-4L^\alpha)}{
\text{no particle of } \widetilde{X}^{K'} \text{ reaches } K' \text{ on } 
\left[ K'^{5/2}, K'^3 \right]
},
\end{align*}
using that $L^3 - L^\alpha \geq K'^3$ and $(L-K')/(\mu' - \mu) = C_\varepsilon L^{2+\alpha} (1+\petito{1}) \leq K'^{5/2}$ for $L$ large enough, with $C_\varepsilon$ a positive constant depending only on $\varepsilon$.
Thus, we get
\begin{align*}
\P_{(L-1-4L^\alpha)}(C_0')
\leq
\P_{(L-1-4L^\alpha)}(R' = 0)
=
1 - \P_{(K'-1)}(R' \geq 1)
\leq
1 - \frac{\E_{(K'-1)}[R']^2}{\E_{(K'-1)}[R'^2]},
\end{align*}
where $R'$ is the number of particles of $\widetilde{X}^{K'}$ that hit $K'$ between times $K'^{5/2}$ and $K'^3$.
We now want to give a lower bound for $\E_{(K'-1)}[R']^2 / \E_{(K'-1)}[R'^2]$.
Using first Proposition \ref{prop estimation esperance de R} with $\theta = 1$, we get
\begin{align*}
\Eci{(K'-1)}{R'}
& \geq
2 \sqrt{2} \pi K' \e^{-\mu K'}
\e^{\mu (K'-1)} \sin \left( \frac{\pi (K'-1)}{K'} \right) (1 + \petito{1})
\underset{L\to\infty}{\longrightarrow} 2 \sqrt{2} \pi^2 \e^{-\mu}.
\end{align*}
Then, using Proposition \ref{prop majoration esperance de R^2} with $R$ denoting the number of particles of $\widetilde{X}^{K'}$ that hit $K'$ between times $0$ and $K'^3$, we get
\begin{align*}
\Eci{(K'-1)}{R'^2}
& \leq
\Eci{(K'-1)}{R^2} \\
& \leq
C \e^{-\mu K'} \e^{\mu (K'-1)} 
\left( 1 + K' \sin \left( \frac{\pi (K'-1)}{K'} \right) \right)
(1 + \petito{1}) \\
& \underset{L\to\infty}{\longrightarrow} C \e^{-\mu} (1+\pi).
\end{align*}
So, for $L$ large enough, we have $\E_{(K'-1)}[R']^2 / \E_{(K'-1)}[R'^2] \geq c > 0$. 
Coming back to (\ref{ef}), we get that, for $L$ large enough,
\[
h(L) 
\leq 
2L \e^{-L^\alpha}
+
(1-c)^L,
\]
and the result follows.
\end{proof}

\section{Upper bound for \texorpdfstring{$v_L$}{vL}} \label{section upper bound}

As in Section \ref{section lower bound}, we fix $0 < \varepsilon < 1/3$ and we consider all processes with drift $-\mu$, but we set here
\[ \mu \coloneqq \sqrt{2 - \frac{\pi^2}{(1+ 4\varepsilon)^2 L^2}}.\]
We want to show that for $L$ large enough, $\P_{(0)}$-a.s.\@ $\lim \max X^L(t) / t \leq C\varepsilon / L^2$, that is $v_L \leq \mu + C\varepsilon / L^2$, and then the upper bound follows by letting $\varepsilon \to 0$.
Moreover, we set for $j\geq 1$, $L_j \coloneqq (1+j\varepsilon)L$.
Thus, $\mu$ is the drift corresponding to $L_4$ according to the results of Subsection \ref{subsection bbs}.

\subsection{Proof of the upper bound}

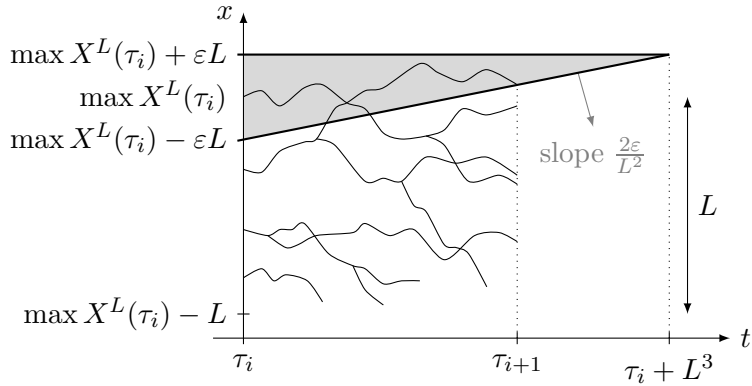
\begin{figure}[ht]
\centering
\begin{tikzpicture}[scale = 0.8]
\fill[gray!30] (0,5.7) -- (7,5.7) -- (0,4.3) -- cycle;

\draw[->,>=latex] (-0.5,1) -- (8,1) node[right]{$t$};
\draw[->,>=latex] (0,0.9) -- (0,6.4) node[left]{$x$};

\draw[->,>=latex,gray] (5.5,4.3 + 1.4/7*5.5) -- (5.75, 4.5) node[below]{slope $\frac{2\varepsilon}{L^2}$};

\draw (0,0.9) node[below]{$\tau_i$};
\draw (7,1.1) -- (7,0.9) node[below]{$\tau_i + L^3$};
\draw (4.5,1.1) -- (4.5,0.9) node[below]{$\tau_{i+1}$};
\draw[dotted] (7,1) -- (7,5.7);
\draw[dotted] (4.5,1) -- (4.5,5.2);

\draw (-0.1,5) node[left]{$\max X^L (\tau_i)$};
\draw (-0.1,5.7) node[left]{$\max X^L (\tau_i) + \varepsilon L$};
\draw (-0.1,4.3) node[left]{$\max X^L (\tau_i) - \varepsilon L$};
\draw (0.1,1.4) -- (-0.1,1.4) node[left]{$\max X^L (\tau_i) - L$};

\draw[thick] (-0.1,5.7) -- (7,5.7);
\draw[thick] (-0.1,4.3-2/70) -- (7,5.7);

\draw[<->,>=latex] (7.3,1.4) -- (7.3,5);
\draw (7.3,3.2) node[right]{$L$};

\draw[smooth] plot coordinates 
{(0,5) (0.3,5.1) (0.6,4.85) (1,5.15) (1.3,5.2) (1.7,4.9)
(2,4.8) (2.3,4.4) (2.6,4.3) (3,4.35) (3.3,4.45) (3.6,4.8) (4,4.7) (4.5,4.85)};
\draw[smooth] plot coordinates 
{(3,4.35) (3.3,4.2) (3.7,3.6) (4.2,3.7) (4.5,3.55)};

\draw[smooth] plot coordinates 
{(0,3.8) (0.3,3.75) (0.6,4.1) (1,4.25) (1.2,4.3) (1.5,4.8) (1.7,4.9)
(2,5.1) (2.3,5.15) (2.9,5.55) (3.2,5.4) (3.6,5.25) (4,5.45) (4.5,5.2)};
\draw[smooth] plot coordinates 
{(1.2,4.3) (1.5,4.2) (2,3.6) (2.3,3.65) (2.6,3.6) (3,4) (3.5,3.9) (4,4) (4.5,3.7)};
\draw[smooth] plot coordinates 
{(2.6,3.6) (2.8,3.2) (3,3.1) (3.3,2.55) (3.5,2.5) (3.7,2.4) (4,5.45-3.6)};

\draw[smooth] plot coordinates 
{(0,2.8) (0.4,2.7) (0.6,2.8) (0.8,2.8) (1,2.9) (1.4,2.4) (1.8,2.2) (2,2.15) (2.3,1.9) (2.6,1.95) (2.9,5.55-3.6)};
\draw[smooth] plot coordinates 
{(1.8,2.2) (2,1.8) (2.3,5.15-3.6)};
\draw[smooth] plot coordinates 
{(0.4,2.7) (0.7,2.5) (0.9,2.55) (1.2,2.7) (1.8,2.75) (2,2.9) (2.4,2.7) (2.8,2.8) (3.3,2.3) (3.8,2.7) (4.2,2.75) (4.5,2.6)};

\draw[smooth] plot coordinates 
{(0,2) (0.3,2.15) (0.5,1.95) (0.7,2.05) (1.05,1.95) (1.3,5.2-3.6)};
\end{tikzpicture}
\caption{Representation of a $L$-BBM between times $\tau_i$ and $\tau_{i+1}$.
By definition, $\tau_{i+1}$ is the time where $t \mapsto \max X^L (t)$ leaves the gray area.
The two thick straight lines that delimit the gray area intersect at time $\tau_i + L^3$.
We are here on the event $A_i^c$ because $t \mapsto \max X^L (t)$ leaves the area from below.}
\label{figure definition tau_i upper bound}
\end{figure}

In this subsection, we prove the upper bound in Theorem \ref{theorem}, by postponing to the next subsections the proof of two propositions.
As in Section \ref{section lower bound}, the $L$-BBM will be studied on time intervals of length at most $L^3$ associated to a sequence of stopping times $(\tau_i)_{i\in\N}$ defined by $\tau_0 \coloneqq 0$ and for each $i \in \N$, 
\[ 
\tau_{i+1} 
\coloneqq 
\inf \left\{ 
t \geq \tau_i : 
\max X^L(t) - \max X^L (\tau_i) \notin
\left( - \varepsilon L + \frac{2 \varepsilon}{L^2} (t-\tau_i), \varepsilon L \right) 
\right\}.
\]
So we have $\tau_{i+1} - \tau_i \leq L^3$ (see Figure \ref{figure definition tau_i upper bound}). 
We also define the event $A_i \coloneqq \{\max X^L (\tau_{i+1}) = \max X^L (\tau_i) + \varepsilon L \}$, as a ``bad'' event: on $A_i$, $\max X^L$ can increase quickly between times $\tau_i$ and $\tau_{i+1}$.
Let $K_n \coloneqq \sum_{i=0}^{n-1} \1_{A_i}$ be the number of ``bad'' events $A_i$ that happen before time $\tau_n$.
On the event $A_i^c$, we have $\max X^L(\tau_{i+1}) - \max X^L(\tau_i) = -\varepsilon L +  2\varepsilon (\tau_{i+1} - \tau_i)/L^2$ and, on the event $A_i$, we have $\max X^L(\tau_{i+1}) - \max X^L(\tau_i) = \varepsilon L$. 
Therefore, we get
\begin{align*}
\max X^L(\tau_n) - \max X^L(0)
& =
\sum_{i=0}^{n-1} \1_{A_i} \varepsilon L
+ \sum_{i=0}^{n-1} \1_{A_i^c} 
		\left( -\varepsilon L + \frac{2 \varepsilon}{L^2}(\tau_{i+1} -\tau_i) \right) \\
& \leq
K_n \varepsilon L - (n-K_n) \varepsilon L + \frac{2 \varepsilon}{L^2} \tau_n
\end{align*}
and so 
\begin{equation} \label{fa}
\frac{\max X^L(\tau_n) - \max X^L(0)}{\tau_n}
\leq
\frac{2 \varepsilon}{L^2} 
+ \frac{n}{\tau_n} \left(- \varepsilon L + 2 \frac{K_n}{n} \varepsilon L \right).
\end{equation}
We now need two propositions to conclude.
The first one shows that there cannot be much more than $(1/2 + 1/\varepsilon L)n$ events $A_i$ happening before time $\tau_n$.
Its proof is postponed to Subsections \ref{subsection comparison} and \ref{subsection proposition K_n}. 
\begin{prop} \label{proposition K_n}
For $L$ large enough, $\P_{(0)}$-a.s.\@ we have
\[ \limsup_{n\to\infty} \frac{K_n}{n} \leq \frac{1}{2} + \frac{1}{\varepsilon L}.\]
\end{prop}
\begin{rem}
The constant $1/2$ appears here because we work in the proof as if the lower barrier that define $\tau_{i+1}$ was horizontal too, so that the population size increase between times $\tau_i$ and $\tau_{i+1}$ on event $A_i^c$ is more or less the inverse of the decrease on event $A_i$.
But, if the lower barrier was horizontal and if there was a vertical barrier at time $L^3$, then we would probably have $\tau_{i+1} = \tau_i + L^3$ most of the time, because fluctuations of the $L$-BBM are believed to be of order $\log L$ (and not $\varepsilon L$) on a time scale of $L^3$.
Thus, with the actual lower barrier, event $A_i^c$ should happen most of the time and $\limsup_{n\to\infty} K_n/n$ should be close to zero.
Moreover, we took $\mu$ greater than the presumed value of $v_L$ and this is also favorable to event $A_i^c$.
\end{rem}
The second proposition gives a lower bound for $\tau_n / n$ as $n \to \infty$ and shows that it is much larger than $L^2$.
Its proof is postponed to Subsection \ref{subsection proposition tau_n}. 
Note that by using some results of Berestycki, Berestycki and Schweinsberg \cite{bbs2014} concerning critical BBM with absorption instead of Corollary \ref{corollary proof tau_n}, we could have $ \liminf_{n\to\infty} \tau_n/n \geq c (\varepsilon L)^3$ for some constant $c>0$ (see Remark \ref{rem tau_n with bbs}).
\begin{prop} \label{proposition tau_n}
There exists $\gamma>0$ such that, for $L$ large enough, $\P_{(0)}$-a.s.\@ we have
\[ \liminf_{n\to\infty} \frac{\tau_n}{n} \geq \frac{L^{2+\gamma}}{6}. \]
\end{prop}
We can now conclude the proof of the upper bound.
Using Proposition \ref{proposition K_n}, we get that $\P_{(0)}$-a.s.\@ $\limsup_{n\to\infty} (- \varepsilon L + 2 \frac{K_n}{n} \varepsilon L ) \leq 2$.
Therefore, it follows from (\ref{fa}) that $\P_{(0)}$-a.s.\@
\begin{align} \label{fb}
\lim_{n\to\infty} \frac{\max X^L(\tau_n)}{\tau_n}
\leq
\frac{2 \varepsilon}{L^2} 
+ 2 \limsup_{n\to\infty} \frac{n}{\tau_n} 
\leq \frac{2 \varepsilon}{L^2} + \frac{12}{L^{2+\gamma}},
\end{align}
for $L$ large enough, applying Proposition \ref{proposition tau_n}.
From Proposition \ref{proposition tau_n}, we get that $\P_{(0)}$-a.s.\@ $\tau_n \to \infty$ as $n\to\infty$, so the left-hand side of (\ref{fb}) is equal to $v_L - \mu$.
Thus, we have
\[
v_L 
\leq 
\sqrt{2} - \frac{\pi^2}{2\sqrt{2} (1+ 4\varepsilon)^2 L^2}
+ \frac{2 \varepsilon}{L^2}
+ \petito{\frac{1}{L^2}}
\]
and, letting $\varepsilon \to 0$, we get the upper bound in Theorem \ref{theorem}.

\subsection{Comparison with the BBM in a strip}
\label{subsection comparison}

We use here results of Subsection \ref{subsection bbs} concerning BBM in a strip to show two lemmas that will be useful for the proof of Proposition \ref{proposition K_n} in the next subsection.
To bound the probability of events $A_i$, we will need to control the size of the $L$-BBM at times $\tau_i$. 
For this, we introduce a functional of the $L$-BBM, analogous to the functional $\widetilde{Z}^{K,K_A}$ of the BBM in a strip (see equation (\ref{equation definition Z})): 
recalling that $L_j = (1 + j\varepsilon) L$, we set, for $i\in \N$,
\begin{align*}
U_i &\coloneqq \max X^L(\tau_i) - L_2 \\
S^L_i &\coloneqq 
\sum_{k=1}^{M^L(\tau_i)} \e^{\mu (X_k^L(\tau_i) - U_i)}
\sin \left( \frac{\pi (X_k^L(\tau_i) - U_i)}{L_4} \right).
\end{align*}
It amounts to shift the population at time $\tau_i$ such that the highest particle is at $L_2$ and then to take the value of the functional $\widetilde{Z}^{L_4}$ associated with the shifted population.
We can now state the two lemmas of this subsection, used later for the proof of Proposition \ref{proposition K_n}.
The first one gives a upper bound for the conditional probability of event $A_i$ given $\sF_{\tau_i}$, in terms of $S^L_i$: if the size of the $L$-BBM at time $\tau_i$ is small enough, then $A_i$ is unlikely.
\begin{lem} \label{lemma upper bound 1} 
There exists $C_\varepsilon > 0$ depending only on $\varepsilon$ such that, for $L$ large enough, for all $i\in\N$,
\[ \Ppsqi{(0)}{A_i}{\sF_{\tau_i}} \leq C_\varepsilon L \e^{-\mu L_3} S^L_i. \]
\end{lem}
Therefore, according to Lemma \ref{lemma upper bound 1}, we need to control $S^L_i$ 
in order to bound the probability of $A_i$.
The second lemma controls the conditional expectation of $S^L_{i+1}$ given $\sF_{\tau_i}$ in terms of $S^L_i$, with a poor bound in the general case but with a much more accurate one on the event $A_i$.
Indeed, even if $A_i$ is a ``bad'' event because it involves a growth of $\max X^L$, it causes at the same time a large decrease of the population size:
when $\max X^L$ grows quickly, more particles are killed by selection.
So each event $A_i$ that happens makes the following events $A_j$ less likely.
\begin{lem} \label{lemma upper bound 2}
We have the following inequalities for $L$ large enough and for all $i\in\N$: 
\begin{align}
\Ecsqi{(0)}{S^L_{i+1}}{\sF_{\tau_i}} 
&\leq 
2 \e^{\mu \varepsilon L} S^L_i \label{inequality S^L 1}, \\
\Ecsqi{(0)}{S^L_{i+1} \1_{A_i}}{\sF_{\tau_i}} 
&\leq 
2 \e^{-\mu \varepsilon L} S^L_i. \label{inequality S^L 2}  
\end{align}
\end{lem}
\begin{proof}[Proof of Lemma \ref{lemma upper bound 1}]
The strategy is as follows: we first come down to the study of $\P_\xi(A_0)$, where $\max \xi = L_2$ and $\min \xi \geq 2\varepsilon L$, and then use Proposition \ref{prop estimation esperance de R} concerning BBM in a strip to bound from above the mean number of particles that hit $L_3$ and thus the probability of event $A_0$.
Applying the strong Markov property at the stopping time $\tau_i$, we get
\begin{equation} \label{gz}
\Ppsqi{(0)}{A_i}{\sF_{\tau_i}} 
=
\Ppi{X^L (\tau_i)}{A_0}
=
\Ppi{X^L (\tau_i) - U_i}{A_0},
\end{equation}
using for the second equality the fact that $\P_\xi (A_0)$ is invariant under shift of the initial configuration $\xi$.
We have $\max X^L (\tau_i) - U_i = L_2$ by definition of $U_i$ and therefore $\min X^L (\tau_i) - U_i \geq 2\varepsilon L$, so we have to bound $\P_\xi(A_0)$ with an initial configuration $\xi$ that satisfies $\max \xi = L_2$ and $\min \xi \geq 2\varepsilon L$.
We fix such a configuration $\xi$.
Then, $\P_\xi$-a.s.\@, for all $t \in [0,\tau_1)$, we have the inclusion $X^L(t) \subset \widetilde{X}^{L_4,L_3}(t)$, because until time $\tau_1$ the killing barrier of $X^L$ stays above 0 (see Figure \ref{figure comparison})
and particles of $X^L$ stay below $L_3$ so no particle of $\widetilde{X}^{L_4,L_3}$ is killed by hitting $L_3$.
Therefore, we get that $\P_\xi$-a.s.\@
\begin{equation} \label{ga}
A_0
\subset
\left\{ 
\text{at least one particle of } \widetilde{X}^{L_4,L_3} \text{ hits } L_3 \text{ on time interval } [0,L^3] 
\right\}.
\end{equation}
Then, 
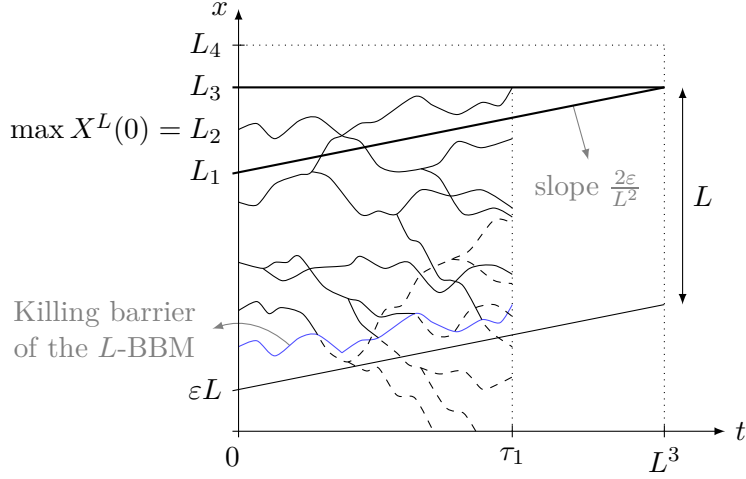
\begin{figure}[t] 
\centering
\begin{tikzpicture}[scale = 0.8]
\node (d) at (1,1.25) {};
\node (f) at (-0.6,1.7) {};
\node[gray,text width=3cm,text centered] (t) at (-2.2,1.7) {Killing barrier of the $L$-BBM};
\draw[->,>=latex,gray] (d) to[bend right] (f.east);
\draw[smooth,blue!70] plot coordinates 
{(0,5-3.6) (0.3,5.1-3.6) (0.6,4.85-3.6) (1,5.15-3.6) (1.3,5.2-3.6) (1.7,4.9-3.6)};
\draw[smooth,blue!70] plot coordinates 
{(1.7,4.9-3.6) (2,5.1-3.6) (2.3,5.15-3.6) (2.9,5.55-3.6) (3.2,5.4-3.6) (3.6,5.25-3.6) (4,5.45-3.6)  (4.3,5.4-3.6) (4.5,5.7-3.6)};

\draw[->,>=latex] (-0.1,0) -- (8,0) node[right]{$t$};
\draw[->,>=latex] (0,-0.1) -- (0,7) node[left]{$x$};

\draw (-0.1,-0.1) node[left,below]{$0$};
\draw (7,0.1) -- (7,-0.1) node[below]{$L^3$};
\draw[dotted] (7,0) -- (7,6.4);

\draw[->,>=latex,gray] (5.5,4.3 + 1.4/7*5.5) -- (5.75, 4.5) node[below]{slope $\frac{2\varepsilon}{L^2}$};

\draw (-0.1,5) node[left]{$\max X^L(0)= L_2$};
\draw (-0.1,5.7) node[left]{$L_3$};
\draw (-0.1,4.3) node[left]{$L_1$};
\draw[dotted] (7,6.4) -- (0,6.4);
\draw (0.1,6.4) -- (-0.1,6.4) node[left]{$L_4$};
\draw (-0.1,0.7) node[left]{$\varepsilon L$};

\draw[thick] (-0.1,5.7) -- (7,5.7);
\draw[thick] (-0.1,4.3-2/70) -- (7,5.7);
\draw (-0.1,0.7-2/70) -- (7,2.1);

\draw (4.5,0.1) -- (4.5,-0.1) node[below]{$\tau_1$};
\draw[dotted] (4.5,0) -- (4.5,5.7);

\draw[<->,>=latex] (7.3,2.1) -- (7.3,5.7);
\draw (7.3,3.9) node[right]{$L$};

\draw[smooth] plot coordinates 
{(0,5) (0.3,5.1) (0.6,4.85) (1,5.15) (1.3,5.2) (1.7,4.9)
(2,4.8) (2.3,4.4) (2.6,4.3) (3,4.35) (3.3,4.45) (3.6,4.8) (4,4.7) (4.5,4.85)};
\draw[smooth] plot coordinates 
{(3,4.35) (3.3,4.2) (3.7,3.6) (4.2,3.7) (4.5,3.55)};

\draw[smooth] plot coordinates 
{(0,3.8) (0.3,3.75) (0.6,4.1) (1,4.25) (1.2,4.3) (1.5,4.8) (1.7,4.9)
(2,5.1) (2.3,5.15) (2.9,5.55) (3.2,5.4) (3.6,5.25) (4,5.45) (4.3,5.4) (4.5,5.7)};
\draw[smooth] plot coordinates 
{(1.2,4.3) (1.5,4.2) (2,3.6) (2.3,3.65) (2.6,3.6) (3,4) (3.5,3.9) (4,4) (4.5,3.7)};
\draw[smooth] plot coordinates 
{(2.6,3.6) (2.8,3.2) (3,3.1) (3.3,2.55) (3.5,2.5) (3.7,2.4) (4,5.45-3.6)};

\draw[smooth] plot coordinates 
{(0,2.8) (0.4,2.7) (0.6,2.8) (0.8,2.8) (1,2.9) (1.4,2.4) (1.8,2.2) (2,2.15) (2.3,1.9) (2.6,1.95) (2.9,5.55-3.6)};
\draw[smooth] plot coordinates 
{(1.8,2.2) (2,1.8) (2.3,5.15-3.6)};
\draw[smooth] plot coordinates 
{(0.4,2.7) (0.7,2.5) (0.9,2.55) (1.2,2.7) (1.8,2.75) (2,2.9) (2.4,2.7) (2.8,2.8) (3.3,2.3) (3.8,2.7) (4.2,2.75) (4.5,2.6)};

\draw[smooth] plot coordinates 
{(0,2) (0.3,2.15) (0.5,1.95) (0.7,2.05) (1.05,1.95) (1.3,5.2-3.6)};

\draw[smooth, dashed] plot coordinates 
{(4+0.05,5.45-3.6-0.09) (4.2,1.6) (4.5,1.4)};

\draw[smooth, dashed] plot coordinates 
{(2.9+0.05,5.55-3.6-0.05) (3.2,1.6) (3.5,1.6) (3.8,2) (4.2,2.1) (4.5,1.9)};

\draw[smooth, dashed] plot coordinates 
{(2.3+0.05,5.15-3.6-0.05) (2.6,1.3) (2.8,1.1) (3,1.05) (3.3,0.95) (3.7,1) (4,0.7) (4.2,0.75) (4.5,0.9)};
\draw[smooth, dashed] plot coordinates 
{(3,1.05) (3.4,0.4) (3.7,0.3) (3.9,0)};

\draw[smooth, dashed] plot coordinates 
{(1.3+0.05,5.2-3.6-0.09) (1.5,1.3) (1.8,1.15) (2.1,1.4) (2.4,2) (2.7,2.1) (3,2.6) (3.5,2.85) (3.8,3) (4.2,3.5) (4.5,3.45)};
\draw[smooth, dashed] plot coordinates 
{(3.5,2.85) (3.8,2.8) (4.2,2.5) (4.5,2.45)};
\draw[smooth, dashed] plot coordinates 
{(1.8,1.15) (2.1,1) (2.3,1.05) (2.6,0.6) (2.9,0.55) (3.2,0)};
\end{tikzpicture}
\caption{Representation of the coupled systems $X^L$ (full line) and $\widetilde{X}^{L_4,L_3}$ (dashed line) between times $0$ and $\tau_1$, starting with the highest particle at $L_2$, on the event $A_0$.
Until time $\tau_1$, the killing barrier of $X^L$ (drawn in blue) stays above the straight line $t\mapsto \varepsilon L + \frac{2\varepsilon}{L^2}t$ so above 0.}
\label{figure comparison}
\end{figure}
our aim is to come down to $\widetilde{X}^{L_3}$ instead of $\widetilde{X}^{L_4,L_3}$, in order to apply Proposition \ref{prop estimation esperance de R}.
For this, let
\[ \mu' \coloneqq \sqrt{2 - \frac{\pi^2}{L_3^2}}\]
be the drift associated to $L_3$ (see equation (\ref{equation mu K})). 
We have $\mu' < \mu$, so the right-hand side of (\ref{ga}) is obviously included in the event
\begin{equation} \label{gb}
\left\{ 
\text{at least one particle of } \widetilde{X}^{L_4,L_3} \text{ hits } 
t \mapsto L_3 +(\mu' - \mu)t \text{ on time interval } [0,L^3]  
\right\}.
\end{equation}
As in the proof of Proposition \ref{prop control P(A) lower bound}, we define now a new process $\bar{X}$ from the standard BBM $X$ with drift $-\mu$ by killing particles that go below $t \mapsto (\mu' - \mu) t$ or above $t \mapsto L_3 + (\mu' - \mu) t$.
Thus, the event in (\ref{gb}) is included in the event
\[
\left\{ 
\text{at least one particle of } \bar{X} \text{ hits } 
t \mapsto L_3 +(\mu' - \mu)t \text{ on time interval } [0,L^3]  
\right\},
\]
because as long as no particle of $\bar{X}$ reaches $t \mapsto L_3 +(\mu' - \mu)t$, the population of $\widetilde{X}^{L_4,L_3}$ is included in the population of $\bar{X}$.
The probability of this last event under $\P_\xi$ is equal to
\[
\Ppi{\xi}
{\text{at least one particle of } \widetilde{X}^{L_3} \text{ hits } 
L_3 \text{ on time interval } [0,L^3] }.
\]
Therefore, coming back to (\ref{ga}), we showed that $\P_\xi(A_0) \leq \E_\xi[R]$, where $R$ is the number of particles of $\widetilde{X}^{L_3}$ that hit $L_3$ between times $0$ and $L^3$, 
and applying Proposition \ref{prop estimation esperance de R} with $\theta = 1 / (1+3\varepsilon)^3$, we get that, as $L\to\infty$, 
\begin{align} \label{gc}
\Eci{\xi}{R}
& \leq
C \frac{\e^{-\mu L_3}}{L_3} \sum_{k=1}^{n} \xi_k \e^{\mu \xi_k}
+
2 \sqrt{2} \pi \theta L_3 \e^{-\mu L_3} 
\left( \sum_{k=1}^{n} \e^{\mu \xi_k} \sin \left( \frac{\pi \xi_k}{L_3} \right) \right)
(1 + \petito{1})
\end{align}
where we wrote $\xi = (\xi_1,\dots,\xi_n)$.
In order to make $S^L_i$ appear, our upper bound has to depend on 
$\sum_{k=1}^{n} \e^{\mu \xi_k} 
\sin(\pi \xi_k/L_4)$.
Recalling that $\xi_k \in [2 \varepsilon L,L_2]$ for each $1 \leq k \leq n$, we get
\[
\frac{\xi_k}{L_2} 
\leq 
\frac{\sin \left( \frac{\pi \xi_k}{L_3} \right)}
	{\sin \left( \frac{\pi L_2}{L_3} \right)}
=
\frac{\sin \left( \frac{\pi \xi_k}{L_3} \right)}
	{\sin \left( \frac{\pi \varepsilon}{(1+3\varepsilon)} \right)}
\]
and it follows that
\begin{equation} \label{gd}
\sum_{k=1}^{n} \xi_k \e^{\mu \xi_k} 
\leq 
\frac{L_2}{\sin \left( \frac{\pi \varepsilon}{(1+3\varepsilon)} \right) } 
\sum_{k=1}^{n} \e^{\mu \xi_k} \sin \left( \frac{\pi \xi_k}{L_3} \right).
\end{equation}
Moreover, using the inequality $\sin(\pi \xi_k /L_3) / \sin(\pi \xi_k /L_4) \leq L_4/L_3$, we have
\begin{equation} \label{ge}
\sum_{k=1}^{n} \e^{\mu \xi_k} \sin \left( \frac{\pi \xi_k}{L_3} \right)
\leq 
\frac{1+4\varepsilon}{1+3\varepsilon} 
\sum_{k=1}^{n} \e^{\mu \xi_k} \sin \left( \frac{\pi \xi_k}{L_4} \right).
\end{equation}
Thus, bringing together (\ref{gc}), (\ref{gd}) and (\ref{ge}), we get that, for $L$ large enough,
\[
\Ppi{\xi}{A_0}
\leq
C_\varepsilon L \e^{-\mu L_3} 
\sum_{k=1}^{n} \e^{\mu \xi_k} \sin \left( \frac{\pi \xi_k}{L_4} \right)
\]
and coming back to (\ref{gz}), the result follows with $\xi = X^L(\tau_i) - U_i$.
\end{proof}
\begin{proof}[Proof of Lemma \ref{lemma upper bound 2}]
The reasoning is similar to the proof of Lemma \ref{lemma upper bound 1}: 
we first come down to the study of $S^L_1$ under $\P_\xi$ where $\max \xi = L_2$, 
then we compare $S^L_1$ with $\widetilde{Z}^{L_4}(\tau_1)$ under $\P_\xi$ in the general case and on the event $A_0$
and finally we use Proposition \ref{prop Z martingale} and the optional stopping theorem to get an upper bound for $\E_\xi [\widetilde{Z}^{L_4}(\tau_1)]$.
Applying the strong Markov property at the stopping time $\tau_i$ and using that the laws of $S^L_1$ and $S^L_1 \1_{A_0}$ are invariant under shift of the initial configuration, we get
\begin{equation} \label{gy}
\Ecsqi{(0)}{S^L_i}{\sF_{\tau_i}} 
=
\Eci{X^L (\tau_i)}{S^L_1}
=
\Eci{X^L (\tau_i) - U_i}{S^L_1}
\end{equation}
and, in the same way, 
$\E_{(0)}[S^L_i \1_{A_i} | \sF_{\tau_i}] = \E_{X^L (\tau_i) - U_i}[S^L_1 \1_{A_0}]$.
Thus, we now have to bound $\E_\xi [S^L_1]$ and $\E_\xi [S^L_1 \1_{A_0}]$, for an initial configuration $\xi$ that satisfies $\max \xi = L_2$.
We fix such a configuration $\xi$.
As in the proof of Lemma \ref{lemma upper bound 1}, we note that $\P_\xi$-a.s.\@, for all $t \in [0,\tau_1]$\footnote{Note that here the inclusion is still true at time $\tau_1$, because a particle of $\widetilde{X}^{L_4}$ that hits $L_3$ is not killed, unlike a particle of $\widetilde{X}^{L_4,L_3}$.}, we have the inclusion $X^L(t) \subset \widetilde{X}^{L_4}(t)$ (see Figure \ref{figure comparison}). 
Thus, we have $\P_\xi$-a.s.\@
\begin{equation} \label{gf}
\sum_{k=1}^{M^L(\tau_1)} \e^{\mu X^L_k(\tau_1)} \sin \left( \frac{\pi X^L_k(\tau_1)}{L_4} \right)
\leq
\widetilde{Z}^{L_4}(\tau_1).
\end{equation}
Therefore, we want to compare $S^L_1$ with the left-hand side of (\ref{gf}), first in the general case and then on the event $A_0$.
Note that $U_1 \in [-\varepsilon L,\varepsilon L]$ and for $1 \leq k \leq M^L(\tau_1)$, $X_k^L(\tau_1) \in [\max X^L(\tau_1) - L,\max X^L(\tau_1)] =  [U_1 + 2 \varepsilon L,U_1 + L_2]$.
On the one hand, if $u \leq 0$, then the function $x \in (0,\pi+u] \mapsto \sin(x-u)/\sin(x)$ is nonincreasing, so if $U_1 \leq 0$, we get the upper bound
\[
\frac{\sin \left( \frac{\pi (X_k^L(\tau_1) - U_1)}{L_4} \right)}
	{\sin \left( \frac{\pi X_k^L(\tau_1)}{L_4} \right)}
\leq
\frac{\sin \left( \frac{\pi 2 \varepsilon L}{L_4} \right)}
	{\sin \left( \frac{\pi (U_1 + 2 \varepsilon L)}{L_4} \right)}
\leq
\frac{\sin \left( \frac{\pi 2 \varepsilon}{1 + 4\varepsilon} \right)}
	{\sin \left( \frac{\pi \varepsilon}{1 + 4\varepsilon} \right)}
\leq 
2.
\]
On the other hand, if $u \geq 0$, then the function $x \in (u,\pi] \mapsto \sin(x-u)/\sin(x)$ is nondecreasing, so if $U_1 \geq 0$, we get the upper bound
\[
\frac{\sin \left( \frac{\pi (X_k^L(\tau_1) - U_1)}{L_4} \right)}
	{\sin \left( \frac{\pi X_k^L(\tau_1)}{L_4} \right)}
\leq
\frac{\sin \left( \frac{\pi L_2}{L_4} \right)}
	{\sin \left( \frac{\pi (U_1 + L_2)}{L_4} \right)}
\leq
\frac{\sin \left( \frac{\pi 2 \varepsilon}{1 + 4\varepsilon} \right)}
	{\sin \left( \frac{\pi \varepsilon}{1 + 4\varepsilon} \right)}
\leq 
2.
\]
It follows that 
\[
S^L_1 
\leq 
2 \e^{- \mu U_1} 
\sum_{k=1}^{M^L(\tau_1)} 
\e^{\mu X^L_k(\tau_1)} \sin \left( \frac{\pi X^L_k(\tau_1)}{L_4} \right).
\]
Coming back to (\ref{gf}) and using that $U_1 \geq -\varepsilon L$ and, on the event $A_0$, $U_1 = \varepsilon L$, we get
\begin{align} \label{gg}
\begin{split}
\Eci{\xi}{S^L_1}
&\leq 
2 \e^{\mu \varepsilon L} \Eci{\xi}{\widetilde{Z}^{L_4}(\tau_1)}, \\
\Eci{\xi}{S^L_1 \1_{A_0}}
&\leq 
2 \e^{-\mu \varepsilon L} \Eci{\xi}{\widetilde{Z}^{L_4}(\tau_1)}.
\end{split}
\end{align}
Note that the bound for $S^L_1 \1_{A_0}$ is better not because the event $A_0$ is unlikely, but because it involves a large increase of $\max X^L$ while most of the particles stay at the same height (because $\mu$ is chosen such that $\widetilde{Z}^{L_4}$ is a martingale) and so all these particles have a much smaller weight in $S^L_1$.
Finally, applying the optional stopping theorem to $(\widetilde{Z}^{L_4}(t))_{t\geq0}$, which is a martingale by Proposition \ref{prop Z martingale}, and to $\tau_1$, which is a bounded stopping time,
we get 
\begin{align} \label{gh}
\Eci{\xi}{\widetilde{Z}^{L_4}(\tau_1)} 
&\leq 
\sum_{k=1}^{n} \e^{\mu \xi_k} \sin \left( \frac{\pi \xi_k}{L_4} \right).
\end{align}
The result follows from (\ref{gy}), (\ref{gg}) and (\ref{gh}) with $\xi = X^L (\tau_i) - U_i$.
\end{proof}

\subsection{Proof of Proposition \ref{proposition K_n}}
\label{subsection proposition K_n}

We prove here Proposition \ref{proposition K_n} that states that $K_n$, the number of events $A_i$ that happen before time $\tau_n$, cannot be much larger than $n/2$.
\begin{proof}[Proof of Proposition \ref{proposition K_n}]
For $1 \leq k \leq n$, we first give a upper bound for $\P_{(0)}(K_n \geq k)$. 
We have
\[
\Ppi{(0)}{K_n \geq k} 
\leq 
\sum_{1 \leq i_1 < \dots < i_k \leq n} \Ppi{(0)}{A_{i_1} \cap \dots \cap A_{i_k}}.
\]
So we fix $1 \leq i_1 < \dots < i_k \leq n$ and deal with $\P_{(0)} (A_{i_1} \cap \dots \cap A_{i_k})$. 
The strategy is to control $S^L_i$ for $1 \leq i \leq i_k$ using Lemma \ref{lemma upper bound 2} 
and then to bound $\P_{(0)}(A_{i_k} | \sF_{\tau_{i_k}})$ with Lemma \ref{lemma upper bound 1} and our control of $S^L_{i_k}$: 
if $k$ is large, then $S^L_{i_k}$ is small and $A_{i_k}$ is unlikely.
First conditioning on $\sF_{\tau_{i_k}}$, using that for all $i\geq0$, $A_i \in \sF_{\tau_{i+1}}$ and then applying Lemma \ref{lemma upper bound 1}, we get
\begin{align} 
\Ppi{(0)}{A_{i_1} \cap \dots \cap A_{i_k}}
& =
\Eci{(0)}{\1_{A_{i_1}} \dotsm \1_{A_{i_{k-1}}} 
	\Ppsqi{(0)}{A_{i_k}}{\sF_{\tau_{i_k}}}} \nonumber \\
& \leq
C_\varepsilon L \e^{-\mu L_3} \Eci{(0)}{\1_{A_{i_1}} \dotsm \1_{A_{i_{k-1}}} S^L_{i_k} }.
	\label{gi}
\end{align}
Then, conditioning successively on all $\sF_i$ for $i$ from $i_{k}-1$ to 0 and applying Lemma \ref{lemma upper bound 2} (we use (\ref{inequality S^L 1}) if $i\notin \{i_1,\dots,i_{k-1}\}$ and (\ref{inequality S^L 2}) otherwise), 
we bound (\ref{gi}) by
\[ 
C_\varepsilon L \e^{-\mu L_3} 
\left( 2 \e^{\mu \varepsilon L} \right)^{i_{k} - (k-1)}
\left( 2 \e^{-\mu \varepsilon L} \right)^{k-1}
\Eci{(0)}{S^L_0}
\]
and it follows that
\begin{equation} \label{gj}
\Ppi{(0)}{K_n \geq k}
\leq
\binom{n}{k}
C_\varepsilon L
\left( 2 \e^{\mu \varepsilon L} \right)^{n - (k-1)}
\left( 2 \e^{-\mu \varepsilon L} \right)^{k-1},
\end{equation} 
using that $\E_{(0)}[S^L_0] = \e^{\mu L_2} \sin(\pi L_2 / L_4) \leq \e^{\mu L_3}$.
We take now
\[ 
k\coloneqq 
\left\lceil n \left(\frac{1}{2} + \frac{1}{\varepsilon L} \right) \right\rceil.
\]
Using $\binom{n}{k} \leq 2^n$, (\ref{gj}) becomes
\begin{align*}
\Ppi{(0)}{K_n \geq k}
& \leq
C_\varepsilon L 
2^{2n}
\exp \left( \mu \varepsilon L (n - 2(k-1)) \right) 
\leq
C L 
2^{2n}
\exp \left(-2 \mu n \right),
\end{align*}
which is summable because $2 < \e^\mu$ for $L$ large enough, 
so the Borel-Cantelli lemma implies that $\P_{(0)}$-almost surely for $n$ large enough we have $K_n \leq k-1$ and the result follows.
\end{proof}

\subsection{Proof of Proposition \ref{proposition tau_n}}
\label{subsection proposition tau_n}

In this section, we fix $0 < \gamma < 1/7$ and we will show that, on the event $A_i^c$, with high probability we have $\tau_{i+1} - \tau_i \geq L^{2+\gamma}$: 
the standard BBM starting at $\max X^L(\tau_i)$ at time $\tau_i$ has with high probability a particle that stays above $t \mapsto \max X^L(\tau_i) - \varepsilon L + 2\varepsilon (t-\tau_i) / L^2$ between times $\tau_i$ and $\tau_i + L^{2+\gamma}$.
Then, by Proposition \ref{proposition K_n}, we know that $\liminf_{n\to\infty} n^{-1} \sum_{i=0}^{n-1} \1_{A_i^c} \geq 1/2 - 1/\epsilon L > 1/3$: thus, more than $n/3$ events $A_i^c$ happen until time $\tau_n$, so we have often enough $\tau_{i+1} - \tau_i \geq L^{2+\gamma}$ and Proposition \ref{proposition tau_n} will follow.
Therefore, we first show a lemma concerning the standard BBM starting with a single particle (we still work with drift $-\mu$).
\begin{lem} \label{lemma proof tau_n}
We define $T \coloneqq L^\gamma + L^{2-5\gamma}$ and the event
\[ 
C \coloneqq 
\left\{ 
\exists k \in \llbracket 1,M(T) \rrbracket : 
X_k(T) \geq - 5 L^\gamma 
\text{ and } 
\forall t \in [0,T], X_{k,T}(t) \geq - \frac{\varepsilon L}{2}
\right\},
\]
where $(X_{k,T}(t))_{0 \leq t \leq T}$ denotes the trajectory between times 0 and $T$ of the particle that is at $X_k(T)$ at time $T$.
Then, for $L$ large enough, we have $\P_{(0)}(C) \geq 1 - 3L\e^{-L^\gamma}$.
\end{lem}
\begin{proof}
The strategy is to use first Lemma \ref{lemma branchement du BBM} in order to get more than $L$ particles after a short time $L^\gamma$ and then Proposition \ref{proposition ABK} to see that each of these particles is likely to stay high between times $L^\gamma$ and $T$.
Applying Lemma \ref{lemma branchement du BBM} (we still have $\mu \leq \sqrt{2}$),
we get
\begin{align} \label{ha}
\Ppi{(0)}{C^c}
& \leq
2 L \e^{-L^\gamma}
+
\Pp{
\left\{ M(L^\gamma) \geq L \right\}
\cap
\left\{ \forall t \in [0,L^\gamma], \min X(t) > -4L^\gamma \right\}
\cap C^c}.
\end{align}
Then, using the branching property at time $L^\gamma$, the probability on the right-hand side of (\ref{ha}) is equal to
\begin{align} \label{hb}
\Eci{(0)}{
\1_{M(L^\gamma) \geq L} 
\1_{\forall t \in [0,L^\gamma], \min X(t) > -4L^\gamma}
\prod_{i=1}^{M(L^\gamma)} 
\Ppi{(X_i(L^\gamma))}{C_1^c} },
\end{align}
for $L$ large enough, where we set
\[
C_1 \coloneqq 
\left\{ 
\exists k \in \llbracket 1,M(L^{2-5\gamma}) \rrbracket : 
X_k(L^{2-5\gamma}) \geq - 5 L^\gamma 
\text{ and } 
\min X_{k,L^{2-5\gamma}}
\geq - \frac{\varepsilon L}{2}
\right\}.
\]
We now want to bound $\P_{(x)}(C_1^c)$ for $x \geq -4L^\gamma$.
First, as $x \mapsto \P_{(x)}(C_1^c)$ is nonincreasing, it is clear that 
$\P_{(x)}(C_1^c) \leq \P_{(-4L^\gamma)}(C_1^c) = \P_{(0)}(C_2^c)$, where we set
\[
C_2 \coloneqq 
\left\{ 
\exists k \in \llbracket 1,M(L^{2-5\gamma}) \rrbracket : 
X_k(L^{2-5\gamma}) \geq - L^\gamma 
\text{ and } 
\min X_{k,L^{2-5\gamma}}
\geq - \frac{\varepsilon L}{2} + 4 L^\gamma
\right\}.
\]
Then, we apply Proposition \ref{proposition ABK} to $\delta = 1/2$: there exist $d>0$, $r>0$ and $t_0 >0$ large enough such that for all $t\geq t_0$, $\P_{(0)}(D_t) \geq 1/2$, where we set for all $t\geq 0$
\begin{align*}
\begin{split}
D_t
& \coloneqq 
\Bigl\{
\exists k \in \llbracket 1, M(t) \rrbracket: 
X_k(t) \geq m(t) - \mu t - d \\
& \hphantom{\coloneqq \Bigl\{} \text{ and } 
\forall s \in [0,t], 
X_{k,t}(s) \geq \frac{s}{t} m(t)  - \mu s
	- r \vee \left(s^{\frac{1}{2}+\gamma} \wedge (t-s)^{\frac{1}{2}+\gamma} \right)
\Bigr\}.
\end{split}
\end{align*}
Note that $D_{L^{2-5\gamma}} \subset C_2$ for $L$ large enough\footnote{On the one hand, with $t \coloneqq L^{2-5\gamma}$ and $0\leq s \leq t$, we have 
$t^{\frac{1}{2}+\gamma} = L^{1 - \frac{\gamma}{2} - 5\gamma^2} = \petito{L}$ 
and $\frac{s}{t} m(t)  - \mu s \geq (\sqrt{2} - \mu)s \geq 0$ 
so $\frac{s}{t} m(t)  - \mu s 	- r \vee (s^{\frac{1}{2}+\gamma} \wedge (t-s)^{\frac{1}{2}+\gamma} ) \geq - \frac{\varepsilon L}{2} + 4 L^\gamma$ for $L$ large enough.
On the other hand, $m(t) - \mu t - d \geq - \frac{3(2-5\gamma)}{2\sqrt{2}} \log L - d \geq - L^\gamma$ for $L$ large enough.},
so we have showed that for $x \geq -4L^\gamma$, $\P_{(x)}(C_1^c) \leq \P_{(0)}(D_{L^{2-5\gamma}}^c) \leq 1/2$ for $L$ large enough such that $L^{2-5\gamma} \geq t_0$ (because $2-5\gamma > 0$).
Thus, we get $(\ref{hb}) \leq 1/2^L$ and, combining with (\ref{ha}), the result follows.
\end{proof}
We now state a corollary that will be used in the proof of Proposition \ref{proposition tau_n}: it says that a standard BBM with drift $-\mu$ starting with a single particle at 0 has with high probability a particle that stays above $t \mapsto - \varepsilon L + 2 \varepsilon t / L^2$ between times 0 and $L^{2+\gamma}$.
\begin{cor} \label{corollary proof tau_n}
We define the event
\[ 
E \coloneqq 
\left\{ 
\exists k \in \llbracket 1,M(L^{2+\gamma}) \rrbracket : 
\forall t \in [0,L^{2+\gamma}], 
X_{k,L^{2+\gamma}}(t) > - \varepsilon L + \frac{2 \varepsilon}{L^2} t 
\right\},
\]
where $(X_{k,L^{2+\gamma}}(t))_{0 \leq t \leq L^{2+\gamma}}$ denotes the trajectory between times 0 and $T$ of the particle that is at $X_k(L^{2+\gamma})$ at time $L^{2+\gamma}$.
Then, $\P_{(0)}(E)$ tends to 1 as $L \to \infty$.
\end{cor}
\begin{proof}
First note that the slope $2 \varepsilon / L^2$ plays only a negligible role on a time period of length $L^{2+\gamma}$ and that we can replace $L^{2+\gamma}$ by $N T$ where $N \coloneqq \lceil L^{6\gamma} \rceil$: 
we have $\bar{E} \subset E$ with
\[ 
\bar{E}
\coloneqq 
\left\{ 
\exists k \in \llbracket 1,M(NT) \rrbracket : 
\forall t \in [0,NT], 
X_{k,NT}(t) \geq - \varepsilon L + L^\gamma
\right\}.
\]
We now want to apply Lemma \ref{lemma proof tau_n} to $N$ consecutive time intervals of length $T$.
Formally, we introduce for $0 \leq j \leq N$ the event
\[ 
\bar{E}_j
\coloneqq 
\left\{ 
\exists k \in \llbracket 1,M(jT) \rrbracket : 
\min X_{k,jT} \geq - \varepsilon L + (5(N-j) + 1) L^\gamma 
\right\}
\]
and we have then, for all $0 \leq j \leq N-1$ and for $L$ large enough such that $(5N +1) L^\gamma \leq \varepsilon L / 2$ (which is possible since $\gamma < 1/7$),
\begin{align} \label{hc}
\begin{split}
\bar{E}_{j+1}
& \supset
\Bigl\{ 
\exists i \in \llbracket 1,M(T) \rrbracket : 
X_i(T) \geq - 5 L^\gamma 
\text{ and } 
\forall t \in [0,T], X_{i,T}(t) \geq - \varepsilon L/2  \text{ and} \\
& \quad\quad
\exists k \in \llbracket 1,M^i(jT) \rrbracket :
X_i(T) + \min X^i_{k,jT} \geq - \varepsilon L + (5(N-j-1) + 1) L^\gamma 
\Bigr\},
\end{split}
\end{align} 
where $X^i$ denotes the BBM emanating from the particle at $X_i(T)$ at time $T$ in the BBM $X$, that is then shifted so that it starts from $0$ at time 0.
According to the branching property, $(X^i, 1 \leq i \leq M^L(T))$ is a family of independent BBM, which is moreover independent of $\sF_T$, so it follows from (\ref{hc}) that 
$\P_{(0)}(\bar{E}_{j+1}) \geq \P_{(0)}(C) \P_{(0)}(\bar{E}_j)$, where $C$ is defined in Lemma \ref{lemma proof tau_n}.
As $E_0 = \Omega$ and $\bar{E}_N = \bar{E}$, we get
\begin{align*}
\P_{(0)}(E) 
&\geq 
\P_{(0)}(C)^N
\geq 
(1 - 2 L \e^{-L^\gamma})^{\lceil L^{4\gamma} \rceil}
=
\exp \left( - 2 L \e^{-L^\gamma} L^{4\gamma} (1 + \petito{1} \right)
\underset{L\to\infty}{\longrightarrow}
1,
\end{align*}
for $L$ large enough, using Lemma \ref{lemma proof tau_n}.
\end{proof}

\begin{rem} \label{rem tau_n with bbs}
Instead of showing Lemma \ref{lemma proof tau_n} and Corollary \ref{corollary proof tau_n}, we could have used Theorem 2 of Berestycki, Berestycki and Schweinsberg \cite{bbs2014}, showing that a BBM with drift $-\sqrt{2}$ starting with a single particle at $x$ has an extinction time close to $(2 \sqrt{2} / 3 \pi^2) x^3$ when $x$ is large enough. 
Thus, Corollary \ref{corollary proof tau_n} is still true with $c (\varepsilon L)^3$ instead of $L^{2 + \gamma}$, for any $c \in (0, 2 \sqrt{2} / 3 \pi^2)$.
But the proof given here is much more elementary and is sufficient for our purpose, so we kept it.
\end{rem}

\begin{proof}[Proof of Proposition \ref{proposition tau_n}]
We introduce for $i\geq 0$ the event $E_i$ defined by ``in the BBM $X$, the particle at $\max X^L (\tau_i)$ at time $\tau_i$ has a descendant at time $\tau_i + L^{2+\gamma}$ whose trajectory between times $\tau_i$ and $\tau_i + L^{2+\gamma}$ stays above $t \mapsto \max X^L (\tau_i) - \varepsilon L + \frac{2 \varepsilon}{L^2} (t-\tau_i)$''.
It is clear that $\P_{(0)}(E_i | \sF_{\tau_i}) = \P_{(0)}(E)$, where $E$ is defined in Corollary \ref{corollary proof tau_n}.
Moreover, we have $A_i^c \cap E_i \subset \{ \tau_{i+1} - \tau_i \geq L^{2+\gamma}\}$: 
on the event $A_i^c$, $\tau_{i+1}$ is the first time after $\tau_i$ when $\max X^L$ goes below $t \mapsto \max X^L (\tau_i) - \varepsilon L + \frac{2 \varepsilon}{L^2} (t-\tau_i)$, 
so on the event $A_i^c \cap E_i$ the descendant at time $\tau_i + L^{2+\gamma}$ in the definition of $E_i$ cannot be killed by selection\footnote{The killing barrier of the $L$-BBM stays below $\max X^L (\tau_i) + \varepsilon L - L$ and, thus, below $t \mapsto \max X^L (\tau_i) - \varepsilon L + \frac{2 \varepsilon}{L^2} (t-\tau_i)$.} 
and so belongs to the $L$-BBM and guarantees that $\tau_{i+1} \geq L^{2+\gamma}+\tau_i$.
Thus, we have
\begin{align} \label{hd}
\Ppsqi{(0)}{A_i^c \cap \{ \tau_{i+1} - \tau_i < L^{2+\gamma}\}}{\sF_{\tau_i}} 
\leq 
\Ppsqi{(0)}{E_i^c}{\sF_{\tau_i}} 
= 
1 - \P_{(0)}(E)
\end{align}
for all $i\geq 0$.
Now the reasoning is as follows: 
by (\ref{hd}) and Corollary \ref{corollary proof tau_n}, on each event $A_i^c$ we have $\tau_{i+1} - \tau_i \geq L^{2+\gamma}$ with high probability 
and, by Proposition \ref{proposition K_n}, we know that more than $n/3$ events $A_i^c$ happen until time $\tau_n$, thus $\tau_n/n$ must be larger than $L^{2+\gamma}/6$ for $n$ large enough.
For $n\in \N^*$, we have
\begin{equation} \label{he}
\Ppi{(0)}{\frac{\tau_n}{n} < \frac{L^{2+\gamma}}{6}}
\leq
\Ppi{(0)}{\frac{\tau_n}{n} < \frac{L^{2+\gamma}}{6} 
	\text{ and } \frac{K_n}{n}\leq\frac{2}{3}}
+ \Ppi{(0)}{\frac{K_n}{n} > \frac{2}{3}}.
\end{equation}
By the proof of Proposition \ref{proposition K_n}, we know that $\P_{(0)} (K_n/n > 2/3)$ is summable in $n$ for $L$ large enough.
Moreover, denoting by $\cP_k(S)$ the set of the subsets with $k$ elements of a set $S$, the first term on the right-hand side of (\ref{he}) is equal to
\begin{align}
& \sum_{k=0}^{\lfloor 2n /3 \rfloor} 
\Ppi{(0)}{\frac{\tau_n}{n} < \frac{L^{2+\gamma}}{6} \text{ and } K_n =k} 
	\nonumber \\
& =
\sum_{k=0}^{\lfloor 2n /3 \rfloor} 
\sum_{I \in \cP_{n-k}(\llbracket 1,n \rrbracket)}
\Ppi{(0)}{ \left\{\frac{\tau_n}{n} < \frac{L^{2+\gamma}}{6} \right\}
\cap \bigcap_{i\in I} A_i^c
\cap \bigcap_{i\notin I} A_i}. \label{hf}
\end{align}
But, on the event $\left\{\tau_n /n < L^{2+\gamma}/6 \right\}$, the events $\{ \tau_{i+1} - \tau_i \geq L^{2+\gamma}\}$ happen for at most $\lfloor n/6\rfloor$ indices $i$. 
Therefore, (\ref{hf}) is bounded by
\begin{align}
& 
\sum_{k=0}^{\lfloor 2n /3 \rfloor} 
\sum_{I \in \cP_{n-k}(\llbracket 1,n \rrbracket)}
\sum_{J \in \cP_{n-k-\lfloor n/6\rfloor}(I)}
\Ppi{(0)}{ 
\bigcap_{i\in J} \left( A_i^c \cap \{ \tau_{i+1} - \tau_i < L^{2+\gamma}\} \right)}
	\nonumber \\
& \leq
\sum_{k=0}^{\lfloor 2n /3 \rfloor} 
\binom{n}{n-k}
\binom{n-k}{n-k-\lfloor n/6\rfloor}
(1-\P_{(0)}(E))^{n-k-\lfloor n/6\rfloor}, \label{hg}
\end{align}
by conditioning successively on $\sF_{\tau_i}$ for all $i\in J$ (in descending order) and using repeatedly (\ref{hd}).
Then, we bound (\ref{hg}) from above by
\begin{align*}
\left( \left\lfloor \frac{2n}{3} \right\rfloor +1 \right) 
2^n
2^n
(1-\P_{(0)}(E))^{n-\lfloor 2n /3 \rfloor-\lfloor n/6\rfloor} 
\leq
n \left(4 (1-\P_{(0)}(E))^{1/6} \right)^n,
\end{align*}
which is summable for $L$ large enough according to Corollary \ref{corollary proof tau_n}.
Coming back to (\ref{he}), 
$\P_{(0)}(\tau_n/n < L^{2+\gamma}/6)$ is summable and the result follows by the Borel-Cantelli lemma.
\end{proof}

\bibliographystyle{abbrv}

\end{document}